\documentclass[a4paper]{article}

\pdfoutput=1

\usepackage{mystyle}
\usepackage{subfiles}

\title{Monoidal Pull-Push I: Cocartesian Fibrations and Categories of Spans}
\author{Angus Rush}

\begin{document}

\maketitle

\tableofcontents

\section{Introduction}
\label{sec:introduction}

\subsection{Introduction}
\label{ssc:introduction}

In topology, one often studies spaces via algebraic invariants associated to them. One such invariant of particular interest is the \emph{$n$th homology} of a space $X$, denoted $H_{n}(X)$. Roughly speaking, $H_{n}(X)$ is the abelian group consisting of formal sums of certain maps from the standard $n$-simplex $\Delta^{n}$ into $X$, considered modulo some equivalence relation.

It turns out to be fruitful to generalize $n$th homology. Formal sums are nothing else but $\Z$-linear combinations, and one can ask what happens when $\Z$ is replaced by some other commutative ring $R$. The associated invariant is then known as \emph{$n$th homology with coefficients in $R$}. One can even allow the ring $R$ to vary along $X$, or allow $R$ to be not only a ring, but an object of some appropriate category $\category{C}$. The structure which keeps track of these changing coefficients is called a $\category{C}$-\emph{local system}.

Classically speaking, when we say \emph{space,} we really mean \emph{topological space,} and $\category{C}$-local systems are modelled by locally-constant sheaves with values in $\category{C}$. In modern homotopy theory, one often takes \emph{space} to mean \emph{$\infty$-groupoid;} thought of in this way, a $\category{C}$-local system on a space $X$ is nothing more than a functor $X \to \category{C}$. This is the point of view that we will take.

With this point of view in mind, let $\category{C}$ be a cocomplete $\infty$-category. For any space $X$, denote by $\LS(\category{C})_{X}$ the $\infty$-category of $\category{C}$-local systems on $X$; that is, $\LS(\category{C})_{X} \simeq \Fun(X, \category{C})$. Recall that for any morphism $f\colon X \to Y$ of spaces, there are several associated functors between $\LS(\category{C})_{X}$ and $\LS(\category{C})_{Y}$. In particular:
\begin{itemize}
  \item The \emph{pullback functor} $f^{*}\colon \LS(\category{C})_{Y} \to \LS(\category{C})_{X}$ pulls back local systems on $Y$ to local systems on $X$, sending a local system $\mathcal{F}\colon Y \to \category{C}$ to the local system $f^{*}\mathcal{F}$ given by the composition
    \begin{equation*}
      \begin{tikzcd}
        X
        \arrow[r, "f"]
        & Y
        \arrow[r, "\mathcal{F}"]
        & \category{C}
      \end{tikzcd}.
    \end{equation*}

  \item The \emph{pushforward functor} $f_{!}\colon \LS(\category{C})_{X} \to \LS(\category{C})_{Y}$ pushes forward local systems via left Kan extension, sending a local system $\mathcal{G}\colon X \to \category{C}$ to the left Kan extension $f_{!}\mathcal{G}$.
    \begin{equation*}
      \begin{tikzcd}
        X
        \arrow[rr, "\mathcal{G}"]
        \arrow[dr, swap, "f"]
        && \category{C}
        \\
        & Y
        \arrow[ur, swap, "f_{!}G"]
      \end{tikzcd}
    \end{equation*}
\end{itemize}

This is a prototype of a common situation: one has a functor (say, of quasicategories) $F\colon \category{C} \to \category{D}$, which sends an object $c \in \category{C}$ to an object $F(c) \in \category{D}$, and a morphism $f\colon c \to c'$ in $\category{C}$ to a morphism $f_{!}\colon F(c) \to F(c')$ in $\category{D}$. One has additionally a `wrong-way' map, which creates from a morphism $f\colon c \to c'$ a morphism $f^{*}\colon F(c') \to F(c)$. In this case, from the data of a diagram in $\category{C}$ of the form
\begin{equation}
  \label{eq:span}
  \begin{tikzcd}
    & c'
    \arrow[dl, swap, "g"]
    \arrow[dr, "f"]
    \\
    c
    && c''
  \end{tikzcd},
\end{equation}
one can produce a morphism $F(c) \to F(c'')$ in $\category{D}$ via the composition
\begin{equation*}
  \begin{tikzcd}
    & F(c')
    \arrow[dr, "f_{!}"]
    \\
    F(c)
    \arrow[ur, "g^{*}"]
    \arrow[rr, swap, "f_{!} \circ g^{*}"]
    && F(c'')
  \end{tikzcd}.
\end{equation*}

The data of \hyperref[eq:span]{Diagram~\ref*{eq:span}} is known as a \emph{span} in $\category{C}$. It is natural to ask whether this construction can be extended to a functor $\Span(\category{C}) \to \category{D}$, where $\Span(\category{C})$ is an $\infty$-category whose objects are the same as the objects of $\category{C}$, and whose morphisms are spans in $\category{C}$. In defining $\Span(\category{C})$, it is necessary to specify a composition law for spans. There is a natural way of doing this: given two spans $X \leftarrow Y \rightarrow X'$ and $X' \leftarrow Y' \rightarrow X''$ in $\category{C}$, we define their composition to be the pullback $X \leftarrow Y \times_{X'} Y' \rightarrow X''$ as below.
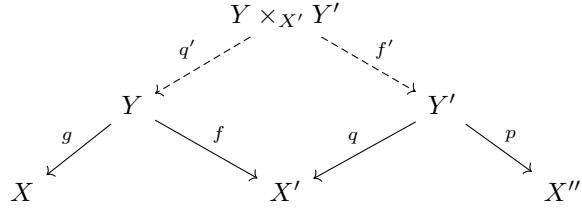
\begin{equation*}
  \begin{tikzcd}
    && Y \times_{X'} Y'
    \arrow[dl, dashed, swap, "q'"]
    \arrow[dr, dashed, "f'"]
    \\
    & Y
    \arrow[dl, swap, "g"]
    \arrow[dr, "f"]
    && Y'
    \arrow[dl, swap, "q"]
    \arrow[dr, "p"]
    \\
    X
    && X'
    && X''
  \end{tikzcd}
\end{equation*}
In order for a functor $\Span(\category{C}) \to \category{D}$ to be well-defined, it will have to respect this composition law on $\Span(\category{C})$. That is, we must have for all such pullback diagrams an equivalence
\begin{equation*}
  (p \circ f')_{!} \circ (g \circ q')^{*} \simeq p_{!} \circ q^{*} \circ f_{!} \circ g^{*}.
\end{equation*}
This is equivalent to the simpler condition that for any pullback square as above we must have an equivalence
\begin{equation*}
  f'_{!} \circ (q')^{*} \simeq q^{*} \circ f_{!}.
\end{equation*}
This is know as the \emph{base change condtition,} or the \emph{Beck--Chevalley condition.}

We will define two $\infty$-categorical models for our category of spans in $\category{C}$: a Segal space $\SPAN(\category{C})$, and a quasicategory $\Span(\category{C})$. Let us for the moment concentrate on the quasicategorical model, although what we say would equally apply to the Segal space model. In defining the quasicategory $\Span(\category{C})$ of spans in $\category{C}$, it is necessary to specify a composition law for spans. There is a natural way of doing this: given two spans $X \leftarrow Y \rightarrow X'$ and $X' \leftarrow Y' \rightarrow X''$ in $\category{C}$, we define their composition to be the pullback $X \leftarrow Y \times_{X'} Y' \rightarrow X''$ as below.
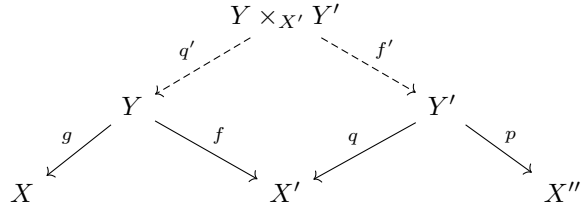
\begin{equation*}
  \begin{tikzcd}
    && Y \times_{X'} Y'
    \arrow[dl, dashed, swap, "q'"]
    \arrow[dr, dashed, "f'"]
    \\
    & Y
    \arrow[dl, swap, "g"]
    \arrow[dr, "f"]
    && Y'
    \arrow[dl, swap, "q"]
    \arrow[dr, "p"]
    \\
    X
    && X'
    && X''
  \end{tikzcd}
\end{equation*}
In order for a functor $\Span(\category{C}) \to \category{D}$ to be well-defined, it will have to respect this composition law on $\Span(\category{C})$. That is, we must have for all such pullback diagrams an equivalence
\begin{equation*}
  (p \circ f')_{!} \circ (g \circ q')^{*} \simeq p_{!} \circ q^{*} \circ f_{!} \circ g^{*}.
\end{equation*}
This is equivalent to the simpler condition that for any pullback square as above we must have an equivalence
\begin{equation*}
  f'_{!} \circ (q')^{*} \simeq q^{*} \circ f_{!}.
\end{equation*}
This is know as the \emph{base change condtition,} or the \emph{Beck--Chevalley condition.}

This is the first part of a two-part paper. In this paper, we to provide a proof of a theorem of Barwick \cite[Thm.~12.2]{spectralmackeyfunctors1}, which provides sufficient conditions for a functor of quasicategories $p\colon \category{C} \to \category{D}$ to yield a cocartesian fibration between $\infty$-categories of spans $\Span(p)\colon \Span(\category{C}) \to \Span(\category{D})$.

In the second part, we will use this result to construct a functor $\hat{r}\colon \Span(\category{S}) \to \ICat$, which sends a space $X$ to the category $\LS(\category{C})_{X}$ of $\category{C}$-local systems on $X$, and a morphism in $\Span(\S)$ represented by a span
\begin{equation*}
  \begin{tikzcd}
    & Y
    \arrow[dl, swap, "g"]
    \arrow[dr, "f"]
    \\
    X
    && X'
  \end{tikzcd}
\end{equation*}
to the functor $f_{!} \circ g^{*}\colon \LS(\category{C})_{X} \to \LS(\category{C})_{X'}$. We further show that our functor is lax monoidal with respect to a monoidal structure on $\Span(\S)$ induced by the cartesian structure on $\S$, and the cartesian structure on $\ICat$.

We begin by providing a proof of the theorem of Barwick mentioned above. Barwick's proof is explicit, constructing the necessary horn fillings by enumerating the necessary simplices, and arguing one-by-one why each filling is possible. This is an impressive feat of simplicial combinatorics, but provided the author little intuition for why the result might be true.

Our proof is more homotopy-theoretic in character, relying on the fact that for any quasicategory $\category{C}$ with pullbacks, the $\infty$-category of spans in $\category{C}$ has a a natural incarnation as a complete Segal space $\SPAN(\category{C})$; the quasicategory $\Span(\category{C})$ is then the `first row' of this complete Segal space. We define a notion of cocartesian fibration between complete Segal spaces, and show that any such cocartesian fibration gives a cocartesian fibration between first rows. The readily available homotopical data in complete Segal spaces allow us to define cocartesian morphisms purely via a condition on 2-simplices, where the combinatorics of horn filling in categories of spans is more manageable.

\subsection{Relation to previous work}
\label{ssc:relation_to_previous_work}

\begin{note}
  This work (MPP I) was originally the author's master's thesis, submitted 21.09.2020. It has come to the author's attention that the same proof of the main theorem (\hyperref[thm:main]{Theorem~\ref*{thm:main}}) was arrived at concurrently as \cite[Thm.~3.1]{2020arXiv201111042H}, with the special case of functors between stable categories appearing already in \cite[Sec.~2.6]{calmes2020hermitian}.
\end{note}

The main result of this paper is a simplified proof of the main theorem of \cite{spectralmackeyfunctors1}, leveraging the theory of cocartesian morphisms in Segal spaces.

Our definition of a cocartesian fibration between Segal spaces is not new, although the form in which it is presented is original. The definition was first written down by De Brito in \cite{2016arXiv160500706B}, by Steimle in \cite{2020arXiv201111042H}, and expanded by Rasekh in \cite{rasekhcartesianfibrations}, and later in \cite{2021arXiv210205192R}. Here, Rasekh proves the existence and equivalence of several model structures for cocartesian fibrations of Segal spaces, including one in terms of marked bisimplicial sets (this did not exist when this paper was written).

In our treatment, also in terms of marked bisimplicial sets, we have tried to balance explicitness of computations, not shying away from combinatorial calculations, with simplicity, not proving results in more generality than we need when we feel it does not lead to greater conceptual clarity. Thus, we have not constructed a model structure of surrounding our definition, and we treat only the case where both the total and base space are fibrant; both of these are treated in full generality in \cite{2021arXiv210205192R}. We feel that the more calculational nature of our treatment leads to a greater interoperability with the existing theory of cocartesian fibrations between quasicategories as found in \cite{highertopostheory}.

\subsection{Outline}
\label{ssc:outline}

This work consists of three sections. In \hyperref[sec:cocartesian_fibrations_between_complete_segal_spaces]{Section~\ref*{sec:cocartesian_fibrations_between_complete_segal_spaces}}, we explore cocartesian fibrations between Segal spaces. After a review in \hyperref[ssc:a_review_of_bisimplicial_sets]{Section~\ref*{ssc:a_review_of_bisimplicial_sets}} of some material in \cite{qcats_vs_segal_spaces}, most notably the box product $- \square -$ and its adjoints, we define marked bisimplicial sets in \hyperref[ssc:marked_bisimplicial_sets]{Section~\ref*{ssc:marked_bisimplicial_sets}}. We then define a marked version of the box functor, and prove some results analogous to the unmarked case. In \hyperref[ssc:simplicial_technology]{Section~\ref*{ssc:simplicial_technology}}, we prove some technical lemmas about simplicial sets with certain restrictions on individual morphisms. The main result is a lemma which allows us to translate a condition on marked bisimplicial sets for a `pointwise condition' involving unmarked bisimplicial sets.

In \hyperref[ssc:cocartesian_fibrations]{Section~\ref*{ssc:cocartesian_fibrations}}, we define the notion of a cocartesian morphism via a condition on left horn filling of $2$-simplices, and define a cocartesian fibration between complete Segal spaces to be a Reedy fibration admitting cocartesian lifts. We show that this implies all higher left horn filling conditions, and use this to show that any cocartesian fibration between complete segal spaces yields a cocartesian fibration (in the sense of quasicategories) between first rows.

In \hyperref[sec:segal_spaces_of_spans]{Section~\ref*{sec:segal_spaces_of_spans}}, we use these results to prove the theorem of Barwick mentioned above, originally published in \cite[Thm.~12.2]{spectralmackeyfunctors1}. The proof leverages the fact, proven in \hyperref[sec:cocartesian_fibrations_between_complete_segal_spaces]{Section~\ref*{sec:cocartesian_fibrations_between_complete_segal_spaces}} that, when checking that a morphism in a Segal space is cocartesian relative to some Reedy fibration, it suffices to check that the lowest lifting problem $\Lambda^{2}_{1} \hookrightarrow \Delta^{2}$ has a contractible space of solutions, which in the setting of categories of spans is combinatorially tractible.

\subsection{Acknowledgements}
\label{ssc:acknowledgements}

We would like to thank his advisor, Tobias Dyckerhoff, for his patience and guidance. We would also like to thank Fernando Abell{\'a}n Garc{\'i}a for many helpful conversations.

\section{Cocartesian fibrations between Segal spaces}
\label{sec:cocartesian_fibrations_between_complete_segal_spaces}

The classical Grothendieck construction allows us to repackage the data of a pseudofunctor $r\colon \category{D} \to \Cat$ into a so-called \emph{cocartesian fibration} $p\colon \category{R} \to \category{D}$. The values of $r$ are encoded in the fibers of $p$ over $\category{D}$. The functoriality of $r$ is encoded by certain morphisms in the total space $\category{R}$, called \emph{cocartesian morphisms.} These morphisms have the following property: \emph{given a commuting triangle in $\category{D}$, and a left horn lying over it in $\category{R}$ whose first edge is cocartesian, there is a unique filling of the horn to a full simplex lying over the triangle.}
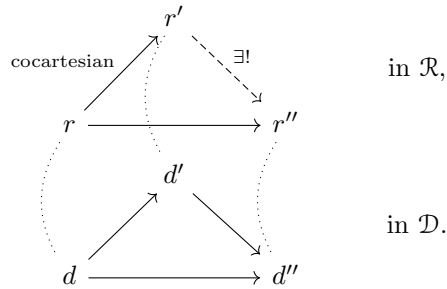
\begin{equation*}
  \label{eq:cocartesian_def}
  \begin{tikzcd}[row sep=tiny]
    & r'
    \arrow[ddd, bend right, dotted, no head]
    \arrow[ddr, dashed, "\exists!"]
    \\
    &&& \text{in $\category{R}$,}
    \\
    r
    \arrow[uur, "\text{cocartesian}"]
    \arrow[rr]
    \arrow[ddd, bend right, dotted, no head]
    && r''
    \arrow[ddd, bend right, dotted, no head]
    \\
    & d'
    \arrow[ddr]
    \\
    &&& \text{in $\category{D}$.}
    \\
    d
    \arrow[uur]
    \arrow[rr]
    && d''
  \end{tikzcd}
\end{equation*}

We can also define cocartesian morphisms in an $\infty$-categorical context. Here, rather than demanding that such fillings should be unique, we should demand that the space of such fillings be contractible. In quasicategories, one does not have easy access to a space of such fillings. The reason for this is that a quasicategory consists of a \emph{set} of 0-simplices, a \emph{set} of 1-simplices, etc., and any homotopical data contained in the quasicategory has to be expressed as a mixing of simplices of various dimensions. For this reason, the most common definition of a cocartesian morphism involves an infinite family of such conditions demanding the existence of horn fillers of ever-increasing dimension.

However one can also model $\infty$-categories as Segal spaces, which one thinks of has having a \emph{space} of 0-simplices, a \emph{space} of 1-simplices, etc. It is then comparably easy to say what one means by the space of diagrams of a specific form. Our definition of a cocartesian edge in a Segal space will then be exactly what we had expected in the $\infty$-categorical case: we will simply replace the \emph{uniqueness} condition on the filler above by a \emph{contractibility} condition.

\subsection{A review of bisimplicial sets}
\label{ssc:a_review_of_bisimplicial_sets}

In this section we review the basic theory of bisimplicial sets as laid out in \cite{qcats_vs_segal_spaces}. This is mainly to fix notation.

Bisimplicial sets can be defined in two equivalent ways:
\begin{itemize}
  \item As functors $\D\op \to \SSet$,

  \item As functors $(\D\op)^{2} \to \Set$.
\end{itemize}
In the former case, we think of a bisimplicial set $X$ as an $\NN$-indexed collection of simplicial sets $X_{n}$; in the latter, we think of a bisimplicial set as an $\NN \times \NN$-indexed collection of sets $X_{mn}$. Both points of view are useful, and we will rely on both of them. For this reason, we fix the following convention: with the second point of view in mind, we imagine a bisimplicial set $X$ as a collection of sets, each located at an integer lattice point of the first quadrant, where the first coordinate increases in the $x$-direction and the second coordinate increases in the $y$-direction.
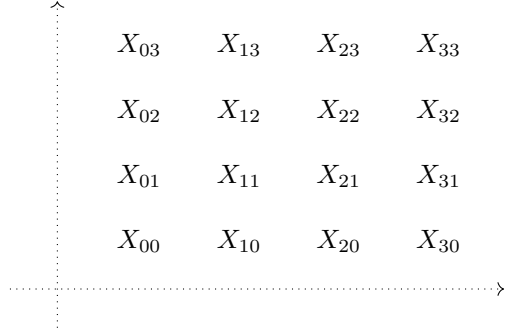
\begin{equation*}
  \begin{tikzcd}[row sep=small, column sep=small]
    & \
    \arrow[dddddd, leftarrow, dotted, at start]
    \\
    && X_{03}
    & X_{13}
    & X_{23}
    & X_{33}
    \\
    && X_{02}
    & X_{12}
    & X_{22}
    & X_{32}
    \\
    && X_{01}
    & X_{11}
    & X_{21}
    & X_{31}
    \\
    && X_{00}
    & X_{10}
    & X_{20}
    & X_{30}
    \\
    \
    \arrow[rrrrrr, dotted, at end]
    &&&&&& \
    \\
    & \
  \end{tikzcd}
\end{equation*}
Thus, the $n$th row of $X$ is the simplicial set $X_{\bullet n}$, and the $m$th column of $X$ is the simplicial set $X_{m \bullet}$. When we think of bisimplicial sets as $\NN$-indexed collections of simplicial sets $X_{m}$, we mean by $X_{m}$ the $m$th \emph{column} of $X$; that is, $X_{m} = X_{m, \bullet}$.

If $X$ is Reedy fibrant (i.e.\ fibrant with respect to the Reedy model structure on $\SSSet$, see \hyperref[ssc:a_review_of_bisimplicial_sets]{Subsection~\ref*{ssc:a_review_of_bisimplicial_sets}}), then each simplicial set $X_{n}$ is a Kan complex. Later, when we are interested in Segal spaces, we will interpret $X_{n}$ as the space of $n$-simplices of $X$.

\begin{definition}
  Define a functor $- \square -\colon \SSet \times \SSet \to \SSSet$ by the formula
  \begin{equation*}
    (X, Y) \mapsto X \square Y,\qquad (X \square Y)_{mn} = X_{m} \times Y_{n}.
  \end{equation*}
  We will call this functor the \defn{box product.}
\end{definition}

\begin{example}
  \label{eg:box_product_and_yoneda}
  By the Yoneda lemma, providing a map $\Delta^{m} \square \Delta^{n} \to X$ is equivalent to providing as an element of $X_{mn}$.
\end{example}

\begin{definition}
  \label{def:left_and_right_divison_functors}
  Let $A$ denote a simplicial set, and $X$ a bisimplicial set.
  \begin{itemize}
    \item Define a simplicial set $A \backslash X$ (pronounced \emph{$A$ under $X$}) level-wise by
      \begin{equation*}
        (A \backslash X)_{n} = \Hom_{\SSet}(A \square \Delta^{n}, X).
      \end{equation*}

    \item Define a simplicial set $X / A$ (pronounced \emph{$X$ over $A$}) level-wise by
      \begin{equation*}
        (X / A)_{n} = \Hom_{\SSet}(\Delta^{n} \square A, X).
      \end{equation*}
  \end{itemize}
\end{definition}

%In \cite{qcats_vs_segal_spaces}, the theory of such bifunctors is explored. The results that we will need are reproduced in \hyperref[sss:divisibility_of_bifunctors]{Appendix~\ref*{sss:divisibility_of_bifunctors}}. However, we note here the following.
Note that by \hyperref[eg:box_product_and_yoneda]{Example~\ref*{eg:box_product_and_yoneda}}, the simplicial set $\Delta^{m} \backslash X$ is the $m$th column of $X$, which we have agreed to call $X_{m}$. Similarly, the simplicial set $X / \Delta^{n}$ is the $n$th row of $X$. In particular $X / \Delta^{0}$ is the zeroth row of $X$. (This, confusingly, is usually called the \emph{first row} of $X$; this terminological inconsistency is somewhat justified by the fact that, unlike the columns, one tends to be interested mainly in the zeroth row.)

\begin{note}
  For $X$ a Reedy-fibrant simplicial set (\hyperref[def:reedy_fibration]{Definition~\ref*{def:reedy_fibration}}), the space $A\backslash X$ is a Kan complex for any simplicial set $A$, and should be thought of as the space of $A$-shaped diagrams in $X$. In particular, $X_{n} = \Delta^{n} \backslash X$ is the space of $n$-simplices in $X$.
\end{note}

We provide here a partial proof of the following result because we will need to refer to it later, when we prove \hyperref[prop:marked_box_adjunctions]{Proposition~\ref*{prop:marked_box_adjunctions}}.
\begin{proposition}
  \label{prop:bijection_exhibiting_box_divisibility}
  The box product is \emph{divisible on the left.} This means that for each simplicial set $A$, there is an adjunction
  \begin{equation*}
    A \square -\colon \SSet \longleftrightarrow \SSSet : A \backslash -.
  \end{equation*}

  Similarly, the box product is \emph{divisible on the right.} This means that for each simplicial set $B$ there is an adjunction
  \begin{equation*}
    - \square B\colon \SSet \longleftrightarrow \SSSet : - / B.
  \end{equation*}
\end{proposition}
\begin{proof}
  We prove divisibility on the left; because the Cartesian product is symmetric, divisibility on the right is identical. We do this by explicitly exhibiting a natural bijection
  \begin{equation*}
    \begin{tikzcd}
      \Hom_{\SSSet}(A \square B, X) \cong \Hom_{\SSet}(B, A \backslash X).
    \end{tikzcd}
  \end{equation*}

  Define a map
  \begin{equation*}
    \Phi\colon \Hom_{\SSSet}(A \square B, X) \to \Hom_{\SSet}(B, A \backslash X)
  \end{equation*}
  by sending a map $f\colon A \square B \to X$ to the map $\tilde{f}\colon B \to A \backslash X$ which sends an $n$-simplex $b \in B_{n}$ to the composition
  \begin{equation*}
    \begin{tikzcd}
      A \square \Delta^{n}
      \arrow[r, "{(\id, b)}"]
      & A \square B
      \arrow[r, "f"]
      & X
    \end{tikzcd}.
  \end{equation*}

  Before we define our map in the other direction, we need an intermediate result. Define a map $\ev\colon A \square (A \backslash X) \to X$ level-wise by taking $(a, \sigma) \in A_{m} \times (A \backslash X)_{n}$ to
  \begin{equation*}
    \sigma_{mn}(a, \id_{\Delta^{n}}) \in X_{mn}.
  \end{equation*}
  Then define a map
  \begin{equation*}
    \Psi\colon \Hom_{\SSet}(B, A \backslash X) \to \Hom_{\SSSet}(A \square B, X)
  \end{equation*}
  sending a map $g\colon B \to A \backslash X$ to the composition
  \begin{equation*}
    \begin{tikzcd}
      A \square B
      \arrow[r, "{(\id, g)}"]
      & A \square (A \backslash X)
      \arrow[r, "\ev"]
      & X
    \end{tikzcd}.
  \end{equation*}

  The maps $\Phi$ and $\Psi$ are mutually inverse, and provide the necessary natural bijection.

  The other bijection is defined analogously, so we only fix notation which we will need later. We will call the mutually inverse maps
  \begin{equation*}
    \Phi'\colon \Hom_{\SSSet}(A \square B, X) \to \Hom_{\SSet}(A, X / B)
  \end{equation*}
  and
  \begin{equation*}
    \Psi'\colon \Hom_{\SSet}(A, X / B) \to \Hom_{\SSSet}(A \square B, X),
  \end{equation*}
  where in defining $\Psi'$ we use a map $\ev'\colon (X / B) \square B \to X$ determined by sending
  \begin{equation*}
    (\phi\colon \Delta^{m} \square B \to X, b \in B_{n}) \mapsto \phi_{mn}(\id_{\Delta^{m}}, b).\qedhere
  \end{equation*}
\end{proof}

\hyperref[prop:bijection_exhibiting_box_divisibility]{Proposition~\ref*{prop:bijection_exhibiting_box_divisibility}}, together with the fact that $\SSet$ and $\SSSet$ are finitely complete and cocomplete, implies all of the results of \hyperref[sss:divisibility_of_bifunctors]{Appendix~\ref*{sss:divisibility_of_bifunctors}} apply to the box product. In the notation found there, we can give a compact formulation of the definition of a Reedy fibration which we will use repeatedly.

\begin{definition}
  \label{def:reedy_fibration}
  Let $f\colon X \to Y$ be a map between bisimplicial sets. The map $f$ is a \defn{Reedy fibration} if either of the following equivalent conditions hold.
  \begin{itemize}
    \item For each monomorphism $u\colon A \to A'$, the map $\langle u \backslash f \rangle$ is a Kan fibration.

    \item For each anodyne map $v\colon B \to B'$, the map $\langle f / v \rangle$ is a trivial Kan fibration.
  \end{itemize}
\end{definition}

(For the notations $\langle u \backslash f \rangle$ and $\langle f / v \rangle$ see \hyperref[sss:divisibility_of_bifunctors]{Appendix~\ref*{sss:divisibility_of_bifunctors}}.) That \hyperref[def:reedy_fibration]{Definition~\ref*{def:reedy_fibration}} is equivalent to the usual definition is shown in \cite[Prop.\ 3.4]{qcats_vs_segal_spaces}. We will also make use of the following fact (\cite[Prop.\ 3.10]{qcats_vs_segal_spaces}).

\begin{theorem}
  \label{thm:inner_fibration_between_quasicategories}
  If $f\colon X \to Y$ is a Reedy fibration between Segal spaces, then for any monomorphism of simplicial sets $v$, the map $\langle f / v \rangle$ is an inner fibration.
\end{theorem}

\begin{corollary}
  \label{cor:reedy_implies_inner}
  If $f\colon X \to Y$ is a Reedy fibration between Segal spaces, then $f / \Delta^{0}\colon X / \Delta^{0} \to Y / \Delta^{0}$ is an inner fibration between quasicategories.
\end{corollary}

\subsection{Marked bisimplicial sets}
\label{ssc:marked_bisimplicial_sets}

In this section, we define a basic theory of marked bisimplicial sets, in analogy to the theory of marked simplicial sets laid out in \cite{highertopostheory}. In the following, one should keep in mind that the case in which we are mostly interested is when our bisimplicial spaces are Segal spaces, and thus that only the first (horizontal) simplicial direction should be thought of as categorical. For this reason, only the first simplicial direction will carry a marking.

\begin{definition}
  A \defn{marked bisimplicial set} $(X, \mathcal{E})$ is a bisimplicial set $X$ together with a distinguished subset $\mathcal{E} \subseteq X_{10}$ containing all degenerate edges, i.e.\ all edges in the image of $s_{0}\colon X_{00} \to X_{10}$. Equivalently, a marked bisimplicial set is bisimplicial set $X$ together with a marking $\mathcal{E}$ on the simplicial set $X / \Delta^{0}$.
\end{definition}

\begin{definition}
  For a marked simplicial set $A$ and an unmarked simplicial set $B$, define a marking on the bisimplicial set $A \square B$ as follows: a simplex $(a, b) \in A_{1} \times B_{0}$ is marked if and only if $a$ is marked in $A$.
\end{definition}

This construction gives us a functor
\begin{equation*}
  - \square -  \colon \SSet^{+} \times \SSet \to \SSSet^{+}.
\end{equation*}

This is potentially ambiguous: we are using the same notation for the marked and unmarked box constructions. However, there is no real chance of confusion: when we write $A \square B$, we mean the marked construction if $A$ is a marked simplicial set and the unmarked construction if $A$ is an unmarked simplicial set.

Our first order of business is to generalize the results of \cite{qcats_vs_segal_spaces} summarized in \hyperref[sss:divisibility_of_bifunctors]{Section~\ref*{sss:divisibility_of_bifunctors}} to the marked case. We will first show that the above functor is divisible on the left and on the right.

\begin{notation}
  For any marked simplicial set $A$, denote the underlying unmarked simplicial set by $\mathring{A}$. Similarly, for any marked bisimplicial set $X$, denote the underlying unmarked bisimplicial set by $\mathring{X}$.
\end{notation}

\begin{definition}
  Let $A$ denote a marked simplicial set, $B$ an unmarked simplicial set, and $X$ a marked bisimplicial set.
  \begin{itemize}
    \item Define an unmarked simplicial set $A \backslash X$ level-wise by
      \begin{equation*}
        (A \backslash X)_{n} = \Hom_{\SSSet^{+}}(A \square \Delta^{n}, X).
      \end{equation*}

    \item Define a marked simplicial set $X / B$ as follows. The underlying simplicial set is the same as $\mathring{X} / B$, and a 1-simplex $\Delta^{1} \to X / B$ is marked if and only if the corresponding map $\Delta^{1} \square B \to \mathring{X}$ of unmarked bisimplicial sets descends to a map of marked bisimplicial sets $(\Delta^{1})^{\sharp} \square B \to X$.
  \end{itemize}
\end{definition}

Again, we are overloading notation, so there is the potential for confusion. However, there is no real ambiguity; the symbol $A \backslash X$ means the marked construction if $A$ and $X$ are marked, and the unmarked construction if $A$ and $X$ are unmarked. We have tried to be clear in stating whether (bi)simplicial sets do or do not carry markings.

\begin{example}
  Recall that we can think of a marked bisimplicial set $X$ as an unmarked bisimplicial set $\mathring{X}$ together with a marking $\mathcal{E}$ on the simplicial set $\mathring{X} / \Delta^{0}$. The marking $\mathcal{E}$ agrees with the marking on $X / \Delta^{0}$.
\end{example}

We will need the following analogs of the $\flat$- and $\sharp$-markings for marked simplicial sets.

\begin{example}
  For any unmarked bisimplicial set $X$, we have the following canonical markings.
  \begin{itemize}
    \item The \emph{sharp marking} $X^{\sharp}$, in which each element of $X_{10}$ is marked.

    \item The \emph{flat marking} $X^{\flat}$, in which only the edges in the image of $s_{0}\colon X_{00} \to X_{10}$ are marked.
  \end{itemize}
\end{example}

\begin{example}
  For each unmarked simplicial set $A$ and marked bisimplicial set $X$, there is an isomorphism
  \begin{equation*}
    A^{\flat} \backslash X \cong A \backslash \mathring{X}.
  \end{equation*}
  Similarly, for any unmarked bisimplicial set $Y$ and marked simplicial set $B$, there is an isomorphism
  \begin{equation*}
    B \backslash Y^{\sharp} \cong \mathring{B} \backslash Y.
  \end{equation*}
\end{example}

The marked constructions above have similar properties to the unmarked constructions from \hyperref[ssc:a_review_of_bisimplicial_sets]{Section~\ref*{ssc:a_review_of_bisimplicial_sets}}. In particular, we have the following.

\begin{proposition}
  \label{prop:marked_box_adjunctions}
  We have the following adjunctions.
  \begin{enumerate}
    \item For each marked simplicial set $A \in \SSet^{+}$ there is an adjunction.
      \begin{equation*}
        A \square -\colon \SSet \longleftrightarrow \SSSet^{+} : A \backslash -
      \end{equation*}

    \item For each unmarked simplicial set $B \in \SSet$ there is an adjunction.
      \begin{equation*}
        - \square B\colon \SSet^{+} \longleftrightarrow \SSSet^{+} : - / B.
      \end{equation*}
  \end{enumerate}
\end{proposition}
\begin{proof}
  We start with the first, fixing a marked simplicial set $A$, an unmarked simplicial set $B$, and a marked bisimplicial set $X$. We have inclusions
  \begin{equation*}
    \Hom_{\SSSet^{+}}(A \square B, X) \overset{i_{0}}{\subseteq} \Hom_{\SSSet}(\mathring{A} \square B, \mathring{X})
  \end{equation*}
  and
  \begin{equation*}
    \Hom_{\SSet}(B, A \backslash X) \overset{i_{1}}{\subseteq} \Hom_{\SSet}(B, \mathring{A} \backslash \mathring{X}).
  \end{equation*}
  We have a natural bijection between the right-hand sides of the above inclusions given by the maps $\Phi$ and $\Psi$ of \hyperref[prop:bijection_exhibiting_box_divisibility]{Proposition~\ref*{prop:bijection_exhibiting_box_divisibility}}. To show that there is a natural bijection between the subsets, it suffices to show that $\Phi$ and $\Psi$ restrict to maps between the subsets.

  To this end, suppose we have a map of marked bisimplicial sets $f\colon A \square B \to X$. The inclusion $i_{0}$ forgets the markings, sending this to the map
  \begin{equation*}
    \mathring{f}\colon \mathring{A} \square B \to \mathring{X}.
  \end{equation*}
  Under $\Phi$, this is taken to a map $\Phi(\mathring{f}) \colon B \to \mathring{A} \backslash \mathring{X}$. We would like to show that $\Phi(\mathring{f})$ factors through $A \backslash X$, giving a map $\tilde{f}\colon B \to A \backslash X$. The map $\Phi(\mathring{f})$ takes an $n$-simplex $b \in B_{n}$ to the composition
  \begin{equation*}
    \begin{tikzcd}
      \mathring{A} \square \Delta^{n}
      \arrow[r, "{(\id, b)}"]
      & \mathring{A} \square B
      \arrow[r, "\mathring{f}"]
      & \mathring{X}
    \end{tikzcd}.
  \end{equation*}
  We need to check that this is an $n$-simplex in $A \backslash X$, and not just $\mathring{A} \backslash \mathring{X}$, i.e.\ that it respects the markings on $A \square \Delta^{n}$ and $X$. That $(\id, b)$ respects the markings on $A \square \Delta^{n}$ and $A \square B$ is clear, and $\mathring{f}$ respects the markings on $A \square B$ and $X$ because $f$ is a map of marked simplicial sets by assumption. Thus $\Phi(\mathring{f})$ restricts to a map $\tilde{f}\colon B \to A \backslash X$.

  Now we show the other direction. Suppose we have a map $g\colon B \to A \backslash X$. The inclusion $i_{1}$ takes this to the composition
  \begin{equation*}
    \begin{tikzcd}
      B
      \arrow[r, "g"]
      & A \backslash X
      \arrow[r, hook]
      & A\flt \backslash X \cong \mathring{A^{\flat}} \backslash \mathring{X},
    \end{tikzcd}
  \end{equation*}
  which we denote by $\mathring{g}$ by mild abuse of notation. Under $\Psi$, this is mapped to the composition
  \begin{equation*}
    \Psi(\mathring{g})\colon
    \begin{tikzcd}
      \mathring{A} \square B
      \arrow[r, "\id \times \mathring{g}"]
      & \mathring{A} \square (\mathring{A} \backslash \mathring{X})
      \arrow[r, "\ev"]
      & \mathring{X}.
    \end{tikzcd}
  \end{equation*}
  We need to check that this respects the markings on $A \square B$ and $X$, i.e.\ that for each marked simplex $a \in A_{1}$ and each $b \in B_{0}$, the element $\Psi(\mathring{g})_{10}(a, b)$ is marked in $X_{10}$. But $\Psi(\mathring{g})_{10}(a, b) = g(b)_{10}(a, \id_{\Delta^{0}})$, which is marked because $g$ lands in $A \backslash X$ by assumption. Thus, $\Psi(\mathring{g})$ descends to a map $\tilde{g}\colon A \square B \to X$.

  Now we show the other bijection. Unlike the unmarked case, because of the asymmetry of the marked box product, this is not precisely the same as what we have just shown. Again we have inclusions
  \begin{equation*}
    \Hom_{\SSSet^{+}}(A \square B, X) \overset{j_{0}}{\subseteq} \Hom_{\SSSet}(\mathring{A} \square B, \mathring{X})
  \end{equation*}
  and
  \begin{equation*}
    \Hom_{\SSet^{+}}(A, B \backslash X) \overset{j_{1}}{\subseteq} \Hom_{\SSet}(\mathring{A}, B \backslash \mathring{X}),
  \end{equation*}
  and a bijection between the right-hand sides given by the maps $\Phi'$ and $\Psi'$ from \hyperref[prop:bijection_exhibiting_box_divisibility]{Proposition~\ref*{prop:bijection_exhibiting_box_divisibility}}. As before, suppose that
  \begin{equation*}
    f\colon A \square B \to X
  \end{equation*}
  is a map of marked bisimplicial sets. Under $j_{0}$, this is sent to $\mathring{f}\colon \mathring{A} \square B \to \mathring{X}$. Then $\Phi'(\mathring{f})\colon A \to X / B$ is defined by sending $\sigma \in A_{n}$ to the composition
  \begin{equation*}
    \begin{tikzcd}
      \Delta^{n} \square B
      \arrow[r, "{(\sigma, \id)}"]
      & \mathring{A} \square B
      \arrow[r, "\mathring{f}"]
      & \mathring{X}
    \end{tikzcd}.
  \end{equation*}
  We need to show that for each marked $a \in A_{1}$, the corresponding map
  \begin{equation*}
    \Phi'(\mathring{f})(a)\colon
    \begin{tikzcd}
      \Delta^{1} \square B
      \arrow[r, "{(a, \id)}"]
      & \mathring{A} \square B
      \arrow[r, "\mathring{f}"]
      & \mathring{X}
    \end{tikzcd}
  \end{equation*}
  descends to a map of marked bisimplicial sets $\tilde{f}\colon (\Delta^{1})^{\sharp} \square B \to X$, and thus corresponds a marked $1$-simplex in to $X / B$. But that the first map has this property is clear because $a$ is marked by assumption, and the map $\mathring{f}$ has this property because $f$ is a map of marked simplicial sets by assumption.

  Now, let $g\colon A \to X / B$ be a map of marked simplicial sets. We need to check that the composition
  \begin{equation*}
    \Psi(\mathring{g})\colon
    \begin{tikzcd}
      \mathring{A} \square B
      \arrow[r, "{(\mathring{g}, \id)}"]
      & (\mathring{X} / B) \square B
      \arrow[r, "\ev'"]
      & \mathring{X}
    \end{tikzcd}
  \end{equation*}
  takes marked edges to marked edges. Let $(a, b) \in A_{1} \times B_{0}$, with $a$ marked. This maps to
  \begin{equation*}
    (a, b) \mapsto (g(a), b) \mapsto g(a)_{10}(\id_{\Delta^{1}}, b) \in X_{10}.
  \end{equation*}
  By definition, $g(a)$ is a map of marked simplicial sets
  \begin{equation*}
    (\Delta^{1})^{\sharp} \square B \to X
  \end{equation*}
  which therefore sends $(\id_{\Delta^{1}}, b)$ to a marked edge in $X$ by assumption.
\end{proof}

This shows that the marked version of the box product $\square$ is, in the language of \cite{qcats_vs_segal_spaces}, \emph{divisible on the left and on the right.} Thus, the results summarized in \hyperref[sss:divisibility_of_bifunctors]{Section~\ref*{sss:divisibility_of_bifunctors}} apply.

\begin{definition}
  \label{def:full_inclusion}
  We will call an inclusion of unmarked simplicial sets $B \hookrightarrow B'$ \emph{full} if it has the following property: an $n$-simplex $\sigma\colon \Delta^{n} \to B'$ factors through $B$ if and only if each vertex of $\sigma$ factors through $B$. That is, any $n$-simplex in $B'$ whose vertices belong to $B$ belongs to $B$.
\end{definition}

\begin{lemma}
  For any marked simplicial set $A$ and Reedy-fibrant marked bisimplicial set $X$, the simplicial set $A \backslash X$ is a Kan complex, and the inclusion $i\colon A \backslash X \hookrightarrow A^{\flat} \backslash X \cong \mathring{A} \backslash \mathring{X}$ is full.
\end{lemma}
\begin{proof}
  We first show that the map $i$ is a full inclusion. The $n$-simplices of $A \backslash X$ are maps of marked simplicial sets $\tilde{\sigma}\colon A \square \Delta^{n} \to X$. A map of underlying bisimplicial sets gives a map of marked bisimplicial sets if and only if it respects the markings, i.e.\ if and only if for each $(a, i) \in A_{1} \times (\Delta^{n})_{0}$ with $a$ marked, $\tilde{\sigma}(a, i)$ is marked in $X$. This is equivalent to demanding that $\sigma|_{\Delta^{\{i\}}}$ belong to $A \backslash X$.

  To show that $A \backslash X$ is a Kan complex, we need to find dashed lifts
  \begin{equation*}
    \begin{tikzcd}
      \Lambda^{n}_{k}
      \arrow[r]
      \arrow[d, hook]
      & A \backslash X
      \\
      \Delta^{n}
      \arrow[ur, dashed]
    \end{tikzcd},
    \qquad n \geq 1,\quad 0 \leq k \leq n.
  \end{equation*}
  For $n = 1$, the horn inclusion is of the form $\Delta^{0} \hookrightarrow \Delta^{1}$, and we can take the lift to be degenerate. For $n \geq 2$, we can augment our diagram as follows.
  \begin{equation*}
    \begin{tikzcd}
      \Lambda^{n}_{k}
      \arrow[r]
      \arrow[d, hook]
      & A \backslash X
      \arrow[r]
      & A^{\flat} \backslash X
      \\
      \Delta^{n}
      \arrow[urr, dashed]
    \end{tikzcd}.
  \end{equation*}
  Since $A^{\flat} \backslash X$ is a Kan complex, we can always find such a dashed lift. The inclusion $\Lambda^{n}_{k} \hookrightarrow \Delta^{n}$ is surjective on vertices, so our lift factors through $A \backslash X$.
\end{proof}

Recall that when thinking of a bisimplicial set $X$ as a simplicial object in $\SSet$, we think of the simplicial set $X_{1}$ as the space of $1$-simplices in $X$. In particular, if $X$ is a Segal space, then $X_{1}$ should be thought of as the space of morphisms in $X$. We should think of morphisms which are in the same path component of $X_{1}$ as equivalent. Therefore, we would like to pay special attention to markings which respect this homotopical structure.
\begin{definition}
  \label{def:respects_path_components}
  Let $(X, \mathcal{E})$ be a marked bisimplicial set. We will say that the marking $\mathcal{E}$ \defn{respects path components} if it has the following property: for any map $\Delta^{1} \to X_{1}$ representing an edge $e \to e'$ between morphisms $e$ and $e'$, the morphism $e$ is marked if and only if the morphism $e'$ is marked.
\end{definition}

\begin{proposition}
  \label{prop:cartesian_marking_respects_path_components}
  Let $f\colon X \to Y$ be a Reedy fibration between marked bisimplicial sets such that the marking on $X$ respects path components, and let $u\colon A \to A'$ be a morphism of marked simplicial sets whose underlying morphism of unmarked simplicial sets is a monomorphism. Then the map $\langle u \backslash f \rangle$ is a Kan fibration.
\end{proposition}
\begin{proof}
  We need to show that for each $n \geq 0$ and $0 \leq k \leq n$ we can solve the lifting problem
  \begin{equation*}
    \begin{tikzcd}
      \Lambda^{n}_{k}
      \arrow[r]
      \arrow[d]
      & A' \backslash X
      \arrow[d, "\langle u \backslash f \rangle"]
      \\
      \Delta^{n}
      \arrow[r]
      \arrow[ur, dashed]
      & A \backslash X \times_{A' \backslash Y} A' \backslash Y
    \end{tikzcd}.
  \end{equation*}
  First assume that $n \geq 2$. We can augment the above square as follows.
  \begin{equation*}
    \begin{tikzcd}
      \Lambda^{n}_{k}
      \arrow[r]
      \arrow[d]
      & A' \backslash X
      \arrow[r]
      \arrow[d]
      & A^{\flat} \backslash X
      \arrow[d]
      \\
      \Delta^{n}
      \arrow[r]
      & A \backslash X \times_{A' \backslash Y} A' \backslash Y
      \arrow[r]
      & A^{\flat} \backslash X \times_{(A')^{\flat} \backslash Y} (A')^{\flat} \backslash X
    \end{tikzcd}.
  \end{equation*}
  Since the map on the right is a Kan fibration, we can solve the outer lifting problem. All the vertices of $\Delta^{n}$ belong to $\Lambda^{n}_{k}$, so a lift of the outside square factors through $A' \backslash X$.

  Now take $n = 1$, $k = 0$, so our horn inclusion is $\Delta^{\{0\}} \hookrightarrow \Delta^{1}$. By \hyperref[prop:equivalent_lifting_problems]{Proposition~\ref*{prop:equivalent_lifting_problems}}, the lifting problem we need to solve is equivalent to
  \begin{equation*}
    \begin{tikzcd}
      A
      \arrow[r]
      \arrow[d]
      & X / \Delta^{1}
      \arrow[d]
      \\
      A'
      \arrow[r]
      \arrow[ur, dashed]
      & X / \Delta^{0} \times_{Y / \Delta^{0}} Y / \Delta^{1}
    \end{tikzcd}.
  \end{equation*}
  Because $f$ is a Reedy fibration, the underlying diagram
  \begin{equation*}
    \begin{tikzcd}
      \mathring{A}
      \arrow[r]
      \arrow[d]
      & \mathring{X} / \Delta^{1}
      \arrow[d]
      \\
      \mathring{A'}
      \arrow[r]
      \arrow[ur, dashed]
      & \mathring{X} / \Delta^{0} \times_{\mathring{Y} / \Delta^{0}} \mathring{Y} / \Delta^{1}
    \end{tikzcd}.
  \end{equation*}
  of unmarked simplicial sets always admits a lift. It therefore suffices to check that any such lift respects the marking on $X$. To see this, consider the following triangle formed by some dashed lift.
  \begin{equation*}
    \begin{tikzcd}[row sep=small, column sep=large]
      & \mathring{X} / \Delta^{1}
      \arrow[dd]
      \\
      \mathring{A'}
      \arrow[ur, dashed]
      \arrow[dr]
      \\
      & \mathring{X} / \Delta^{0}
    \end{tikzcd}
  \end{equation*}
  Let $a \in A'_{1}$ be a marked 1-simplex, and consider the diagram
  \begin{equation*}
    \begin{tikzcd}[row sep=small, column sep=large]
      \Delta^{1} \square \Delta^{0}
      \arrow[r, "{(a, \id)}"]
      \arrow[dd]
      & \mathring{A}' \square \Delta^{0}
      \arrow[dd]
      \arrow[dr, "\gamma"]
      \\
      && \mathring{X}
      \\
      \Delta^{1} \square \Delta^{1}
      \arrow[r, "{(a, \id)}"]
      & \mathring{A}' \square \Delta^{1}
      \arrow[ur, dashed, swap, "\ell"]
    \end{tikzcd},
  \end{equation*}
  where the triangle on the right is the adjunct to the triangle above. In order to check that the dashed lift respects the marking on $X$, we have to show that for each $(a, b) \in (A' \square \Delta^{1})_{10} = A'_{1} \times \{0, 1\}$ with $a$ marked, the element $\ell(a, b) \in X_{10}$ is marked. Because the map $\gamma$ comes from a map of marked simplicial sets, the commutativity of the triangle guarantees this for $b = 0$. The map $\Delta^{1} \square \Delta^{1} \to \mathring{X}$ gives us a 1-simplex $\Delta^{1} \to X_{1}$ representing a 1-simplex $\ell(a, 0) \to \ell(a, 1)$, which implies by that $\ell(a, 1)$ is also marked because each marking respects path components.

  The case $n = 1$, $k = 1$ is analogous.
\end{proof}

\subsection{Simplicial technology}
\label{ssc:simplicial_technology}

In the next section, we will need to work in several different cases with simplicial subsets $A \subseteq \Delta^{n}$ with certain conditions placed on the edge $\Delta^{\{0, 1\}}$. In this section we prove some technical results in this direction. The main result in this section is \hyperref[lemma:check_marked_trivial_fibration_pointwise]{Lemma~\ref*{lemma:check_marked_trivial_fibration_pointwise}}.

For the remainder of this section, fix $n \geq 2$.

\begin{definition}
  \label{def:pullback_of_diagrams}
  Let $X$ be an unmarked bisimplicial set, and let $e \in X_{10}$. For any simplicial subset $A \subseteq \Delta^{n}$ such that $\Delta^{\{0, 1\}} \subseteq A$, we will use the notation
  \begin{equation*}
    (A \backslash X)^{e} = A \backslash X \times_{\Delta^{\{0, 1\}} \backslash X} \{e\}.
  \end{equation*}
\end{definition}

The simplicial set $(A \backslash X)^{e}$ should be thought of as the space of $A$-shaped diagrams in $X$ with the edge $\Delta^{\{0, 1\}}$ fixed. The $m$-simplices of the simplicial set $(A \backslash X)^{e}$ are maps $A \square \Delta^{m} \to X$ such that the pullback
\begin{equation*}
  \begin{tikzcd}
    \Delta^{\{0, 1\}} \square \Delta^{m}
    \arrow[r]
    & A \square \Delta^{m}
    \arrow[r]
    & X
  \end{tikzcd}
\end{equation*}
factors through the map $\Delta^{\{0, 1\}} \square \Delta^{0} \to X$ corresponding to the element $e \in X_{10}$ under the Yoneda embedding.

Comparing simplices level-wise, it is easy to see the following.
\begin{lemma}
  \label{lemma:unmarked_pullback}
  The square
  \begin{equation*}
    \begin{tikzcd}
      (\Delta^{n} \backslash X)^{e}
      \arrow[r]
      \arrow[d]
      & \Delta^{n} \backslash X
      \arrow[d]
      \\
      (\Delta^{n} \backslash Y)^{f(e)} \times_{(A \backslash Y)^{f(e)}} (A \backslash X)^{e}
      \arrow[r]
      & \Delta^{n} \backslash Y \times_{A\backslash Y} A \backslash X
    \end{tikzcd}
  \end{equation*}
  is a (strict) pullback.
\end{lemma}

\begin{definition}
  For any simplicial subset $A \subseteq \Delta^{n}$ containing $\Delta^{\{0, 1\}}$, denote the marking on $A$ where the only marked nondegenerate edge is $\Delta^{\{0, 1\}}$ by $\mathcal{L}$, and the corresponding marked simplicial set by $A^{\mathcal{L}}$.
\end{definition}

Again, comparing simplices level-wise shows the following.
\begin{lemma}
  \label{lemma:marked_pullback}
  Let $f\colon X \to Y$ be a map of marked bisimplicial sets, and let $A \subseteq \Delta^{n}$ be a simplicial subset with $\Delta^{\{0, 1\}} \subseteq A$. Then for any marked edge $e \in X_{10}$, the square
  \begin{equation*}
    \begin{tikzcd}
      (\Delta^{n} \backslash \mathring{X})^{e}
      \arrow[r]
      \arrow[d]
      & (\Delta^{n})^{\mathcal{L}} \backslash X
      \arrow[d]
      \\
      (\Delta^{n} \backslash \mathring{Y})^{f(e)} \times_{(A\backslash \mathring{Y})^{f(e)}} (A \backslash \mathring{X})^{e}
      \arrow[r]
      & (\Delta^{n})^{\mathcal{L}} \backslash Y \times_{A^{\mathcal{L}}\backslash Y} A^{\mathcal{L}} \backslash X
    \end{tikzcd}
  \end{equation*}
  is a (strict) pullback.
\end{lemma}

\begin{lemma}
  \label{lemma:check_marked_trivial_fibration_pointwise}
  Let $f\colon X \to Y$ be a Reedy fibration between marked bisimplicial sets, and let $i\colon A \subseteq \Delta^{n}$ be a simplicial subset containing $\Delta^{\{0, 1\}}$. The following are equivalent:
  \begin{enumerate}
    \item The map
      \begin{equation*}
        \langle i^{\mathcal{L}} \backslash f \rangle\colon (\Delta^{n})^{\mathcal{L}} \backslash X \to (\Delta^{n})^{\mathcal{L}} \backslash Y \times_{A^{\mathcal{L}} \backslash Y} A^{\mathcal{L}} \backslash X
      \end{equation*}
      is a trivial fibration.

    \item For each marked $e \in X_{10}$, the map
      \begin{equation*}
        p_{e}\colon (\Delta^{n} \backslash \mathring{X})^{e} \to (\Delta^{n} \backslash \mathring{Y})^{f(e)} \times_{(A \backslash \mathring{Y})^{f(e)}} (A \backslash \mathring{X})^{e}
      \end{equation*}
      is a trivial fibration.
  \end{enumerate}
\end{lemma}
\begin{proof}
  Suppose the first holds. Then \hyperref[lemma:marked_pullback]{Lemma~\ref*{lemma:marked_pullback}} implies the second.

  Next, suppose that the second holds. By \hyperref[prop:cartesian_marking_respects_path_components]{Proposition~\ref*{prop:cartesian_marking_respects_path_components}}, the map $\langle i^{\mathcal{L}} \backslash f \rangle$ is a Kan fibration, so it is a trivial Kan fibration if and only if its fibers are contractible. Consider any map
  \begin{equation*}
    \gamma\colon \Delta^{0} \to (\Delta^{n})^{\mathcal{L}} \backslash Y \times_{A^{\mathcal{L}} \backslash Y} A^{\mathcal{L}} \backslash X.
  \end{equation*}
  This gives us in particular a map $\Delta^{0} \to A^{\mathcal{L}} \backslash X$, which is adjunct to a map $A^{\mathcal{L}} \square \Delta^{0} \to X$. The pullback
  \begin{equation*}
    \begin{tikzcd}
      (\Delta^{\{0, 1\}})^{\sharp} \square \Delta^{0}
      \arrow[r]
      & A^{\mathcal{L}} \square \Delta^{0}
      \arrow[r]
      & X
    \end{tikzcd}
  \end{equation*}
  gives us a marked morphism $e \in X_{10}$. The bottom composition in the below diagram is thus a factorization of $\gamma$, in which the left-hand square is a pullback.
  \begin{equation*}
    \begin{tikzcd}
      F
      \arrow[r]
      \arrow[d]
      & (\Delta^{n} \backslash \mathring{X})^{e}
      \arrow[r]
      \arrow[d, "p_{e}"]
      & (\Delta^{n})^{\mathcal{L}} \backslash X
      \arrow[d]
      \\
      \Delta^{0}
      \arrow[r]
      & (\Delta^{n} \backslash \mathring{Y})^{f(e)} \times_{(A \backslash \mathring{Y})^{f(e)}} (A \backslash \mathring{X})^{e}
      \arrow[r]
      & (\Delta^{n})^{\mathcal{L}} \backslash Y \times_{A^{\mathcal{L}} \backslash Y} A^{\mathcal{L}} \backslash X
    \end{tikzcd}
  \end{equation*}
  The right-hand square is a pullback by \hyperref[lemma:marked_pullback]{Lemma~\ref*{lemma:marked_pullback}}. Since by assumption $p_{e}$ is a trivial fibration, $F$ is contractible. But by the pasting lemma, $F$ is the fiber of $\langle i^{\mathcal{L}} \backslash f \rangle$ over $\gamma$. Thus, the fibers of $\langle i^{\mathcal{L}} \backslash f \rangle$ are contractible, so $\langle i^{\mathcal{L}} \backslash f \rangle$ is a trivial Kan fibration.
\end{proof}

\subsection{Cocartesian fibrations}
\label{ssc:cocartesian_fibrations}

Let $\pi\colon \category{C} \to \category{D}$ be an inner fibration between quasicategories. There are several equivalent ways of defining when a morphism in $\category{C}$ is $\pi$-cocartesian. For our purposes, the following will be the most useful: a morphism $e \in \category{C}_{1}$ is $\pi$-cocartesian if and only if, for all $n \geq 2$, a dashed lift in the below diagram exists. 
\begin{equation*}
  \begin{tikzcd}
    \Delta^{\{0, 1\}}
    \arrow[d, hook]
    \arrow[dr, "e"]
    \\
    \Lambda^{n}_{0}
    \arrow[r]
    \arrow[d, hook]
    & \category{C}
    \arrow[d, "\pi"]
    \\
    \Delta^{n}
    \arrow[r]
    \arrow[ur, dashed]
    & \category{D}
  \end{tikzcd}
\end{equation*}
We would like to find an analogous definition for a $p$-cocartesian morphism where $p\colon C \to D$ is a Reedy fibration between Segal spaces. Our definition should have the property that if a morphism in $C_{10}$ is $p$-cocartesian in the sense of Segal spaces, then it is $p / \Delta^{0}$-cocartesian in the sense of quasicategories.

To this end, replace $\pi$ in the above diagram by $p / \Delta^{0}$. Passing to the adjoint lifting problem, we see that the existence of the above lift is equivalent to demanding that the map
\begin{equation}
  p_{e}\colon (\Delta^{n} \backslash C)^{e} \to (\Lambda^{n}_{0} \backslash C)^{e} \times_{(\Lambda^{n}_{0} \backslash D)^{p(e)}} (\Delta^{n} \backslash D)^{p(e)}
\end{equation}
be surjective on vertices. One natural avenue of generalization of the concept of a cocartesian morphism to Segal spaces would be to upgrade the condition of surjectivity on vertices to an analogous, homotopy-invariant condition which implies it. One such condition is that $p_{e}$ be a trivial fibration. Indeed, this is the defintion we will use. However, this turns out to be equivalent to demand something superficially weaker.

\begin{definition}
  \label{def:cocartesian_morphism}
  Let $f\colon X \to Y$ be a Reedy fibration between Segal spaces. A morphism $e \in X_{10}$ is \defn{$f$-cocartesian} if the square
  \begin{equation*}
    \begin{tikzcd}
      (\Delta^{2} \backslash X)^{e}
      \arrow[r]
      \arrow[d]
      & (\Lambda^{2}_{0} \backslash X)^{e}
      \arrow[d]
      \\
      (\Delta^{2} \backslash Y)^{f(e)}
      \arrow[r]
      & (\Lambda^{2}_{0} \backslash Y)^{f(e)}
    \end{tikzcd}
  \end{equation*}
  is homotopy pullback.
\end{definition}

\begin{example}
  \label{eg:id_and_equiv_are_cocartesian}
  Identity morphisms are $f$-cocartesian. This is because for $e = \id$, the horizontal morphisms in \hyperref[def:cocartesian_morphism]{Definition~\ref*{def:cocartesian_morphism}} are equivalences. More generally, by \cite[Lemma\ 11.6]{rezk2001model}, homotopy equivalences are $f$-cocartesian.
\end{example}

Cocartesian morphisms automatically respect path components in the following sense.

\begin{proposition}
  \label{prop:cocartesian_morphisms_respect_path_components}
  Let $f\colon X \to Y$ be a Reedy fibration between Segal spaces, and let $\alpha\colon \Delta^{1} \to X_{1}$ be a map representing a path $e \to e'$ between morphisms $e$ and $e'$ in $X_{1}$. Then $e$ is $f$-cocartesian if and only if $e'$ is $f$-cocartesian.
\end{proposition}
\begin{proof}
  Since $X$ and $Y$ are Reedy fibrant and $f$ is a Reedy fibration, it suffices to show that the map
  \begin{equation}
    \label{eq:unprimed_map}
    (\Delta^{2} \backslash X)^{e} \to
    (\Lambda^{2}_{0} \backslash X)^{e}
    \times_{(\Lambda^{2}_{0} \backslash Y)^{f(e)}}
    (\Delta^{2} \backslash Y)^{f(e)}
  \end{equation}
  is a weak equivalence if and only if the map
  \begin{equation}
    \label{eq:primed_map}
    (\Delta^{2} \backslash X)^{e'} \to
    (\Lambda^{2}_{0} \backslash X)^{e'}
    \times_{(\Lambda^{2}_{0} \backslash Y)^{f(e')}}
    (\Delta^{2} \backslash Y)^{f(e')}
  \end{equation}
  is a weak equivalence.

  Consider the diagram
  \begin{equation*}
    \begin{tikzcd}
      P
      \arrow[r]
      \arrow[d]
      & \Delta^{2} \backslash X
      \arrow[d]
      \\
      Q
      \arrow[r]
      \arrow[d]
      & \Lambda^{2}_{0} \backslash X \times_{\Lambda^{2}_{0} \backslash Y} \Delta^{2} \backslash Y
      \arrow[d]
      \\
      \Delta^{1}
      \arrow[r, "\alpha"]
      & \Delta^{\{0, 1\}} \backslash X
    \end{tikzcd},
  \end{equation*}
  where both squares are strict pullback. The maps on the right-hand side are Kan fibrations by Reedy fibrancy and the fact that $f$ is a Reedy fibration, so we get in particular a diagram
  \begin{equation*}
    \begin{tikzcd}[column sep=small]
      P
      \arrow[rr, "\phi"]
      \arrow[dr]
      && Q
      \arrow[dl]
      \\
      & \Delta^{1}
    \end{tikzcd},
  \end{equation*}
  where both downward-facing maps are Kan fibrations. Note that the component of the map $\phi$  over $\Delta^{0}$ is the map from \hyperref[eq:unprimed_map]{Equation~\ref*{eq:unprimed_map}}, and component over $\Delta^{1}$ is the map from \hyperref[eq:primed_map]{Equation~\ref*{eq:primed_map}}. Kan fibrations are in particular left fibrations, so under the Grothendieck construction this corresponds to a diagram
  %of the form
  %\begin{equation*}
  %  \begin{tikzcd}
  %    \mathfrak{C}[\Delta^{1}]
  %    \arrow[bend left ]{r}[name=U, label=above:$\scriptstyle\mathrm{St}_{\Delta^{1}} P$]{}
  %    \arrow[bend right]{r}[name=D, label=below:$\scriptstyle\mathrm{St}_{\Delta^{1}} Q$]{}
  %    & \SSet
  %    \arrow[shorten <=5pt, from=U, to=D, Rightarrow]
  %  \end{tikzcd},
  %\end{equation*}
  %givin
  given by the following homotopy-commutative square in $\SSet$ (with the Kan model structre), in which the rightward-pointing maps are weak equivalences because our maps $P \to \Delta^{1}$ and $Q \to \Delta^{1}$ were Kan fibrations.
  \begin{equation*}
    \begin{tikzcd}
      (\Delta^{2} \backslash X)^{e}
      \arrow[r, "\simeq"]
      \arrow[d]
      & (\Delta^{2} \backslash X)^{e'}
      \arrow[d]
      \\
      (\Lambda^{2}_{0} \backslash X)^{e} \times_{(\Lambda^{2}_{0} \backslash Y)^{f(e)}} (\Delta^{2} \backslash Y)^{f(e)}
      \arrow[r, "\simeq"]
      & (\Lambda^{2}_{0} \backslash X)^{e'} \times_{(\Lambda^{2}_{0} \backslash Y)^{f(e')}} (\Delta^{2} \backslash Y)^{f(e')}
    \end{tikzcd}
  \end{equation*}
  By the $2/3$ property for weak equivalences, the map on the left is a weak equivalence if and only if the map on the right is a weak equivalence, which is what we wanted to show.
\end{proof}

\begin{definition}
  \label{def:morphism_generated_marking}
  Let $f\colon X \to Y$ be a Reedy fibration between Segal spaces, and let $e \in X_{10}$ be a cocartesian morphism. Denote by $\mathcal{E}$ the smallest marking which contains $e$ and all degenerate morphisms in $X$, and which respects path components. More explicitly, the marking $\mathcal{E}$ contains:
  \begin{itemize}
    \item The morphism $e$;

    \item Each identity morphism; and

    \item Any morphism connected to $e$ or any identity morphism by a path.
  \end{itemize}
\end{definition}

Our next order of business is to show that our definition of cocartesian morphisms in terms of lifting with respect to the morphism $\Lambda^{2}_{0} \hookrightarrow \Delta^{2}$ implies lifting with respect to $\Lambda^{n}_{0} \to \Delta^{n}$ for all $n \geq 2$.

\begin{definition}
  Define the following simplicial subsets of $\Delta^{n}$.
  \begin{itemize}
    \item For $n \geq 1$, denote by $I_{n}$ the \defn{spine} of $\Delta^{n}$, i.e.\ the simplicial subset
      \begin{equation*}
        \Delta^{\{0, 1\}} \amalg_{\Delta^{\{1\}}} \Delta^{\{1, 2\}} \amalg_{\Delta^{\{2\}}} \cdots \amalg_{\Delta^{\{n-1\}}} \Delta^{\{n-1, n\}} \subseteq \Delta^{n}.
      \end{equation*}

    \item For $n \geq 2$, denote by $L_{n}$ the simplicial subset
      \begin{equation*}
        L_{n} = \Delta^{\{0, 1\}} \amalg_{\Delta^{\{0\}}} \overbrace{\Delta^{\{0, 2\}} \amalg_{\Delta^{\{2\}}} \Delta^{\{2, 3\}} \amalg_{\Delta^{\{3\}}}\cdots \amalg_{\Delta^{\{n-1\}}} \Delta^{\{n-1, n\}}}^{I_{\{0, \hat{1}, 2, \ldots, n\}}} \subseteq \Delta^{n}.
      \end{equation*}
      That is, $L_{n}$ is the union of $\Delta^{\{0, 1\}}$ with the spine of $d_{1}\Delta^{n}$. We will call $L_{n}$ the \defn{left spine} of $\Delta^{n}$.
  \end{itemize}
\end{definition}

Note that $L_{2} \cong \Lambda^{2}_{0}$.

%\begin{lemma}
%  For any $n \geq 2$, the map
%  \begin{equation*}
%    \Delta^{\{0, \ldots, n-1\}} \amalg_{\Delta^{\{n-1\}}} \Delta^{\{n-1, n\}} \hookrightarrow \Delta^{n}
%  \end{equation*}
%  is inner anodyne.
%\end{lemma}
%\begin{proof}
%  This map can be written as the starred smash product of $\emptyset \hookrightarrow \Delta^{n-2}$, which is an inclusion, and $\Delta^{\{0\}} \hookrightarrow \Delta^{1}$, which is left anodyne. The result follows from \cite[Lemma 2.1.2.3]{highertopostheory}.
%\end{proof}

\begin{proposition}
  \label{prop:only_lowest_lifting_condition_is_necessary}
  Let $f\colon X \to Y$ be a Reedy fibration between Segal spaces, and let $e \in X_{10}$ be an $f$-cocartesian morphism. Then the square
  \begin{equation*}
    \begin{tikzcd}
      (\Delta^{n} \backslash X)^{e}
      \arrow[r]
      \arrow[d]
      & (L_{n} \backslash X)^{e}
      \arrow[d]
      \\
      (\Delta^{n} \backslash Y)^{f(e)}
      \arrow[r]
      & (L_{n} \backslash Y)^{f(e)}
    \end{tikzcd}
  \end{equation*}
  is homotopy pullback for all $n \geq 2$.
\end{proposition}
\begin{proof}
  We have the case $n = 2$ because $e$ is $f$-cocartesian. Assume the result is true up to $n - 1$. Then the square
  \begin{equation*}
    \begin{tikzcd}
      (\Delta^{n-1} \backslash X)^{e} \times_{\Delta^{\{n-1\}} \backslash X} \Delta^{\{n-1, n\}} \backslash X
      \arrow[r]
      \arrow[d]
      & (L_{n-1} \backslash X)^{e} \times_{\Delta^{\{n-1\}} \backslash X} \Delta^{\{n-1, n\}} \backslash X
      \arrow[d]
      \\
      (\Delta^{n-1} \backslash Y)^{f(e)} \times_{\Delta^{\{n-1\}} \backslash Y} \Delta^{\{n-1, n\}} \backslash Y
      \arrow[r]
      & (L_{n-1} \backslash Y)^{f(e)} \times_{\Delta^{\{n-1\}} \backslash Y} \Delta^{\{n-1, n\}} \backslash Y
    \end{tikzcd}
  \end{equation*}
  is homotopy pullback since each component is homotopy pullback. But this square is equivalent to
  \begin{equation*}
    \begin{tikzcd}
      (\Delta^{n} \backslash X)^{e}
      \arrow[r]
      \arrow[d]
      & (L_{n} \backslash X)^{e}
      \arrow[d]
      \\
      (\Delta^{n} \backslash Y)^{f(e)}
      \arrow[r]
      & (L_{n} \backslash Y)^{f(e)}
    \end{tikzcd}:
  \end{equation*}
  The left-hand equivalences come from the Segal condition, and the right-hand equivalences come from the definition of $L_{n}$.
\end{proof}

\begin{corollary}
  \label{cor:marked_left_spine_gives_triv_fib}
  Let $f\colon X \to Y$ be a Reedy fibration between Segal spaces, and let $e \in X_{10}$ be a $f$-cocartesian morphism. Then for all $n \geq 2$, the map
  \begin{equation*}
    (\Delta^{n})^{\mathcal{L}} \backslash X^{\mathcal{E}} \to L_{n}^{\mathcal{L}} \backslash X^{\mathcal{E}} \times_{(\Delta^{n})^{\mathcal{L}} \backslash Y\shp} L_{n}^{\mathcal{L}} \backslash Y\shp
  \end{equation*}
  is a trivial Kan fibration.
\end{corollary}
\begin{proof}
  Each edge $e' \in \mathcal{E}$ is $f$-cocartesian: the morphism $e$ is $f$-cocartesian by assumption, each degenerate edge is $f$-cocartesian by \hyperref[eg:id_and_equiv_are_cocartesian]{Example~\ref*{eg:id_and_equiv_are_cocartesian}}, and any morphism in the path component of an $f$-cocartesian morphism is $f$-cocartesian by \hyperref[prop:cocartesian_morphisms_respect_path_components]{Proposition~\ref*{prop:cocartesian_morphisms_respect_path_components}}.

  Therefore, for any $e' \in \mathcal{E}$, the map
  \begin{equation*}
    (\Delta^{n} \backslash X)^{e'} \to (L_{n} \backslash X)^{e'} \times_{(\Delta^{n} \backslash Y)^{f(e')}} (L_{n} \backslash Y)^{f(e')}
  \end{equation*}
  is a weak equivalence by \hyperref[prop:only_lowest_lifting_condition_is_necessary]{Proposition~\ref*{prop:only_lowest_lifting_condition_is_necessary}}, and it is a Kan fibration by \hyperref[lemma:unmarked_pullback]{Lemma~\ref*{lemma:unmarked_pullback}}. The result follows from \hyperref[lemma:check_marked_trivial_fibration_pointwise]{Lemma~\ref*{lemma:check_marked_trivial_fibration_pointwise}}.
\end{proof}

For any simplicial set $A$, define a marked simplicial set $(\Delta^{1} \star A, \mathcal{L'})$ where the only nondegenerate simplex belonging to $\mathcal{L}'$ is $\Delta^{1}$. This is a slight generalization of the $\mathcal{L}$-marking.

\begin{lemma}
  \label{lemma:starred_smash_with_mono}
  Let $A \hookrightarrow B$ be a monomorphism of simplicial sets, and suppose that $B$ is $n$-skeletal (and therefore that $A$ is $n$-skeletal). Then the map
  \begin{equation*}
    \begin{tikzcd}
      (\Delta^{\{0\}} \star B)^{\flat} \coprod_{(\Delta^{\{0\}} \star A)^{\flat}} (\Delta^{1} \star A)^{\mathcal{L}'} \hookrightarrow (\Delta^{1} \star B)^{\mathcal{L}'}
    \end{tikzcd}
  \end{equation*}
  is in the saturated hull of the morphisms
  \begin{equation*}
    (\Lambda^{k}_{0})^{\mathcal{L}} \hookrightarrow (\Delta^{k})^{\mathcal{L}},\qquad 2 \leq k \leq n+2.
  \end{equation*}
\end{lemma}
\begin{proof}
  It suffices to show this for $A \hookrightarrow B = \partial \Delta^{m} \hookrightarrow \Delta^{m}$ for $0 \leq m \leq n$. In this case the necessary map is of the form
  \begin{equation*}
    (\Lambda^{m+2}_{0})^{\mathcal{L}} \hookrightarrow (\Delta^{m+2})^{\mathcal{L}}.
  \end{equation*}
\end{proof}

\begin{definition}
  \label{def:right_cancellation_property}
  We will say a collection of morphisms $\mathcal{A} \subset \mathrm{Mor}(\SSet^{+})$ has the \defn{right cancellation property} if for all $u$, $v \in \mathrm{Mor}(\SSet^{+})$,
  \begin{equation*}
    u \in \mathcal{A},\quad vu \in \mathcal{A} \quad \implies \quad v \in A.
  \end{equation*}
\end{definition}

\begin{lemma}
  \label{lemma:saturated_hull_of_left_spine_inclusions}
  Let $\mathcal{A}$ be a saturated set of morphisms of $\SSet^{+}$ all of whose underlying morphisms are monomorphisms, and which has the right cancellation property. Further suppose that $\mathcal{A}$ contains the following classes of morphisms.
  \begin{enumerate}
    \item Maps $(A)^{\flat} \hookrightarrow (B)^{\flat}$, where $A \to B$ is inner anodyne.

    \item Left spine inclusions $(L_{n})^{\mathcal{L}} \hookrightarrow (\Delta^{n})^{\mathcal{L}}$, $n \geq 2$.
  \end{enumerate}

  Then $\mathcal{A}$ contains left horn inclusions $(\Lambda^{n}_{0})^{\mathcal{L}} \hookrightarrow (\Delta^{n})^{\mathcal{L}}$, $n \geq 2$.
\end{lemma}
\begin{proof}
  For $n = 2$, there is nothing to check: we have an isomorphism $(L_{2})^{\mathcal{L}} \cong (\Lambda^{2}_{0})^{\mathcal{L}}$.

  We proceed by induction. Suppose we have shown that all horn inclusions $(\Lambda^{k}_{0})^{\mathcal{L}} \hookrightarrow (\Delta^{k})^{\mathcal{L}}$ belong to $\mathcal{A}$ for $2 \leq k < n$. From now on on we will suppress the marking $(-)^{\mathcal{L}}$. All simplicial subsets of $\Delta^{n}$ below will have $\Delta^{\{0, 1\}}$ marked if they contain it.

  Consider the factorization
  \begin{equation*}
    \begin{tikzcd}
      L_{n}
      \arrow[r, "u_{n}"]
      \arrow[rr, bend right, swap, "v_{n} \circ u_{n}"]
      & \Lambda^{n}_{0}
      \arrow[r, "v_{n}"]
      & \Delta^{n}
    \end{tikzcd}.
  \end{equation*}
  The morphism $v_{n} \circ u_{n}$ belongs to $\mathcal{A}$ by assumption, so in order to show that $v_{n}$ belongs to $\mathcal{A}$, it suffices by right cancellation to show that $u_{n}$ belongs to $\mathcal{A}$. Consider the factorization
  \begin{equation*}
    \begin{tikzcd}
      L_{n}
      \arrow[r, "w'_{n}"]
      \arrow[rr, bend right, swap, "u_{n}"]
      & L_{n} \cup d_{1} \Delta^{n}
      \arrow[r, "w_{n}"]
      & \Lambda^{n}_{0}
    \end{tikzcd}.
  \end{equation*}
  The map $w'_{n}$ is a pushout along the spine inclusion $I_{\{0, \hat{1}, 2, \ldots, n\}} \hookrightarrow d_{1}\Delta^{n}$, and hence is inner anodyne. Hence, we need only show that $w_{n}$ belongs to $\mathcal{A}$. Let
  \begin{equation*}
    Q = d_{2} \Delta^{n} \cup \cdots \cup d_{n} \Delta^{n},
  \end{equation*}
  and consider the following pushout diagram.
  \begin{equation*}
    \begin{tikzcd}
      (L_{n} \cup d_{1}\Delta^{n}) \cap Q
      \arrow[r, hook]
      \arrow[d, hook]
      & Q
      \arrow[d, hook]
      \\
      L_{n} \cup d_{1} \Delta^{n}
      \arrow[r, hook]
      & L_{n} \cup d_{1}\Delta^{n} \cup Q
    \end{tikzcd}
  \end{equation*}
  Since $L_{n} \cup d_{1} \Delta^{n} \cup Q \cong \Lambda^{n}_{0}$, the bottom map is $w_{n}$, so it suffices to show that the top map belongs to $\mathcal{A}$. But this is isomorphic to
  \begin{equation*}
    \begin{tikzcd}
      (\Delta^{\{0, 1\}} \star \emptyset) \coprod_{(\Delta^{\{0\}} \star \emptyset)} (\Delta^{\{0\}} \star \partial \Delta^{\{2, 3, \ldots, n\}}) \hookrightarrow \Delta^{\{0, 1\}} \star \partial \Delta^{\{2, 3, \ldots, n\}}.
    \end{tikzcd}
  \end{equation*}
  The simplicial set $\partial \Delta^{\{2, \ldots, n\}}$ is $(n-3)$-skeletal, so this map belongs to $\mathcal{A}$ by \hyperref[lemma:starred_smash_with_mono]{Lemma~\ref*{lemma:starred_smash_with_mono}}.
\end{proof}

For each $n \geq 2$, denote by $h^{n}$ the $\mathcal{L}$-marked inclusion
\begin{equation*}
  h^{n}\colon (\Lambda^{n}_{0})^{\mathcal{L}} \hookrightarrow (\Delta^{n})^{\mathcal{L}}.
\end{equation*}

\begin{proposition}
  \label{prop:segal_cocartesian_morphisms_are_quasicategory_cocartesian}
  Let $f\colon X \to Y$ be a Reedy fibration of Segal spaces, and let $e \in X_{10}$ be an $f$-cocartesian morphism. Then for all $n \geq 2$, the map
  \begin{equation*}
    \langle h^{n} \backslash f^{\mathcal{E}} \rangle\colon (\Delta^{n})^{\mathcal{L}} \backslash X^{\mathcal{E}} \to (\Lambda^{n}_{0})^{\mathcal{L}} \backslash X^{\mathcal{E}} \times_{(\Lambda^{n}_{0})^{\mathcal{L}} \backslash Y^{\sharp}} (\Delta^{n})^{\mathcal{L}} \backslash Y^{\sharp}
  \end{equation*}
  is a trivial fibration of simplicial sets.
\end{proposition}
\begin{proof}
  Consider the set
  \begin{equation*}
    S =
    \left\{
      \substack{
        u\colon A \to \text{ B morphism of} \\
        \text{marked simplicial sets} \\
        \text{such that $\mathring{u}$ is mono}
      }
      \ \bigg| \
      \langle u \backslash f^{\mathcal{E}} \rangle \text{ weak homotopy equivalence}
    \right\}.
  \end{equation*}
  We claim that this set has the right cancellation property (\hyperref[def:right_cancellation_property]{Definition~\ref*{def:right_cancellation_property}}). The set of morphisms of marked simplicial sets whose underlying morphisms are monic clearly has the right-cancellation property. To show that $S$ does, let $u\colon A \to B$ and $v\colon B \to C$ be such morphisms and consider the following diagram.
  \begin{equation*}
    \begin{tikzcd}[column sep=large]
      C \backslash X
      \arrow[r, "\langle vu \backslash f^{\mathcal{E}} \rangle"]
      \arrow[d, swap, "\langle v \backslash f^{\mathcal{E}} \rangle"]
      & A \backslash X \times_{A \backslash Y} C \backslash Y
      \\
      B \backslash X \times_{B \backslash Y} C \backslash Y
      \arrow[r, "\langle u \backslash f^{\mathcal{E}} \rangle \times_{\id} \id"]
      & \left( A \backslash X \times_{A \backslash Y} B \backslash Y \right) \times_{B \backslash Y} C \backslash Y
      \arrow[u, swap, "\simeq"]
    \end{tikzcd}
  \end{equation*}
  If $\langle u \backslash f^{\mathcal{E}} \rangle$ is a weak equivalence, then the bottom morphism is a weak equivalence. The right-hand morphism is a weak equivalence because it is an isomorphism, so if $\langle vu \backslash f^{\mathcal{E}} \rangle$ is a weak equivalence, then the $\langle v \backslash f^{\mathcal{E}} \rangle$ is a weak equivalence by $2 / 3$.

  By \hyperref[prop:cartesian_marking_respects_path_components]{Proposition~\ref*{prop:cartesian_marking_respects_path_components}}, a map $u$ belonging to $S$ automatically has the property that $\langle u \backslash f^{\mathcal{E}} \rangle$ a Kan fibration, hence is a trivial Kan fibration. Thus, we can equivalently say that $u \in S$ if and only if $u$ has the left-lifting property with respect to all maps of the form $\langle X / v \rangle$, where $v$ is a cofibration of simplicial sets. Since the set of all monomorphisms is saturated, $S$ is saturated.

  The set $S$ contains all flat-marked inner anodyne morphisms because $f$ is a Reedy fibration. \hyperref[cor:marked_left_spine_gives_triv_fib]{Corollary~\ref*{cor:marked_left_spine_gives_triv_fib}} tells us that $S$ contains all left spine inclusions $(L_{n})^{\mathcal{L}} \hookrightarrow (\Delta^{n})^{\mathcal{L}}$, $n \geq 2$. Thus, by \hyperref[lemma:saturated_hull_of_left_spine_inclusions]{Lemma~\ref*{lemma:saturated_hull_of_left_spine_inclusions}}, $S$ contains all $\mathcal{L}$-marked left horn inclusions.
\end{proof}

\begin{corollary}
  \label{cor:pointwise_left_horn_condition}
  Let $f\colon X \to Y$ be a Reedy fibration between Segal spaces, and let $e \in X_{10}$ be an $f$-cocartesian edge. Then the map
  \begin{equation*}
    (\Delta^{n} \backslash X)^{e} \to (\Lambda^{n}_{0} \backslash X)^{e} \times_{(\Lambda^{n}_{0} \backslash Y)^{f(e)}} (\Delta^{n} \backslash Y)^{f(e)}
  \end{equation*}
  is a trivial fibration.
\end{corollary}
\begin{proof}
  Apply \hyperref[lemma:check_marked_trivial_fibration_pointwise]{Lemma~\ref*{lemma:check_marked_trivial_fibration_pointwise}} to \hyperref[prop:segal_cocartesian_morphisms_are_quasicategory_cocartesian]{Proposition~\ref*{prop:segal_cocartesian_morphisms_are_quasicategory_cocartesian}}.
\end{proof}

\begin{corollary}
  \label{cor:cocartesian_morphisms_are_first_row_cocartesian}
  Let $f\colon X \to Y$ be a Reedy fibration between Segal spaces, and let $e \in X_{10}$ be an $f$-cocartesian morphism. Then $e$ is $f / \Delta^{0}$-cocartesian.
\end{corollary}
\begin{proof}
  By \hyperref[cor:pointwise_left_horn_condition]{Corollary~\ref*{cor:pointwise_left_horn_condition}}, for all $n \geq 2$ the map
  \begin{equation*}
    (\Delta^{n} \backslash X)^{e} \to (\Lambda^{n}_{0} \backslash X)^{e} \times_{(\Lambda^{n}_{0} \backslash Y)^{f(e)}} (\Delta^{n} \backslash Y)^{f(e)}
  \end{equation*}
  is a trivial fibration. Thus, it certainly has the right-lifting property with respect to $\emptyset \hookrightarrow \Delta^{0}$. Passing to the adjoint lifting problem, we find that this is equivalent to the existence of a dashed lift in the diagram
  \begin{equation*}
    \begin{tikzcd}
      \Delta^{\{0, 1\}}
      \arrow[d, hook]
      \arrow[dr, "e"]
      \\
      \Lambda^{n}_{0}
      \arrow[r]
      \arrow[d, hook]
      & X / \Delta^{0}
      \arrow[d]
      \\
      \Delta^{n}
      \arrow[r]
      \arrow[ur, dashed]
      & Y / \Delta^{0}
    \end{tikzcd},
  \end{equation*}
  which tells us that $e$ is $f / \Delta^{0}$-cocartesian.
\end{proof}

\begin{definition}
  \label{def:cocartesian_fibration_between_complete_segal_spaces}
  Let $f\colon X \to Y$ be a Reedy fibration between Segal spaces. We will say that $f$ is a \defn{cocartesian fibration} if each morphism in $Y$ has an $f$-cocartesian lift in $X$. More explicitly, we demand that, for each edge $e\colon y \to y'$ in $Y$ and each vertex $x \in X$ such that $f(x) = y$, there exists an $f$-cocartesian morphism $\tilde{e}\colon x \to x'$ such that $f(\tilde{e}) = e$.
\end{definition}

\begin{corollary}
  \label{cor:cocart_fib_between_css_gives_cocart_fib_of_quasicats}
  Let $f\colon X \to Y$ be a cocartesian fibration of Segal spaces. Then the map
  \begin{equation*}
    f/\Delta^{0}\colon X/\Delta^{0} \to Y/\Delta^{0}
  \end{equation*}
  is a cocartesian fibration of quasicategories, and if a morphism in $X_{1}$ is $f$-cocartesian, then it is $f / \Delta^{0}$-cocartesian.
\end{corollary}
\begin{proof}
  By \hyperref[thm:inner_fibration_between_quasicategories]{Theorem~\ref*{thm:inner_fibration_between_quasicategories}}, the map $f/\Delta^{0}$ is an inner fibration between quasicategories. By assumption, every morphism in $Y$ has a $f$-cocartesian lift, and these lifts are $f / \Delta^{0}$-cocartesian by  \hyperref[cor:cocartesian_morphisms_are_first_row_cocartesian]{Corollary~\ref*{cor:cocartesian_morphisms_are_first_row_cocartesian}}.
\end{proof}

%\subfile{unconstrained_spans.tex}

\section{Segal spaces of spans}
\label{sec:segal_spaces_of_spans}

\subsection{Basic definitions}

We recall the basic definitions of Segal spaces of spans. Note that spans as we will define them form not only a Segal space, but a \emph{complete} Segal space; for the most part, this will not concern us. For more information, we direct the reader to \cite{spectralmackeyfunctors1}. This section is intended to be a summary of the relevant results of loc.\ cit.

The objects of our study will be $\infty$-categories whose morphisms are spans in some quasicategory $\category{C}$, i.e.\ diagrams in $\category{C}$ of the form
\begin{equation*}
  \begin{tikzcd}
    & y
    \arrow[dl, swap, "\phi"]
    \arrow[dr, "\psi"]
    \\
    x
    && x'
  \end{tikzcd}.
\end{equation*}
We will want to be able to place certain conditions on the legs $\phi$ and $\psi$ of our spans. For example, we may want to restrict our attention to spans such that $\phi$ is an equivalence. More precisely, we pick out two subcategories of $\category{C}$ to which the respective legs of our spans must belong.

\begin{definition}
  A \defn{triple} of categories is a triple $\triple{C}$, where $\category{C}$ is a quasicategory where $\category{C}_{\dagger}$ and $\category{C}^{\dagger}$ are subcategories, each of which contain all equivalences.
\end{definition}

We will use the following terminology and notation for the morphsism in our subcategories.
\begin{itemize}
  \item We denote the morphisms in $\category{C}_{\dagger}$ with tails (as in $x \ing y$), and call them \emph{ingressive}.

  \item We denote the morphisms in $\category{C}^{\dagger}$ with two heads (as in $x \eg y$), and call them \emph{egressive}.
\end{itemize}

The egressive morphisms will correspond to the backwards (i.e.\ leftwards) facing legs of our spans, and the ingressive morphisms will correspond to the forwards (i.e.\ rightwards) facing legs. We provide the following diagram as a summary of this terminology.
\begin{equation*}
  \begin{tikzcd}
    & Y
    \arrow[dl, two heads, swap, "\substack{\text{egressive} \\ \category{C}\updag}"]
    \arrow[dr, rightarrowtail, "\substack{\text{ingressive} \\ \category{C}\downdag}"]
    \\
    X
    && X'
  \end{tikzcd}
\end{equation*}

In order for a triple of categories to be able to support a category of spans, it will have to satisfy certain properties. Following Barwick, we will call such triples \emph{adequate.}

\begin{definition}
  \label{def:adequate_triple}
  A triple $(\category{C}, \category{C}_{\dagger}, \category{C}^{\dagger})$ is said to be \defn{adequate} if it has the following properties:
  \begin{enumerate}
    \item For any ingressive morphism $f \in \category{C}_{\dagger}$ and any egressive morphism $g \in \category{C}^{\dagger}$, there exists a pullback square
      \begin{equation*}
        \begin{tikzcd}
          y'
          \arrow[r]
          \arrow[d]
          & x'
          \arrow[d, two heads, "g"]
          \\
          y
          \arrow[r, rightarrowtail, "f"]
          & x
        \end{tikzcd}.
      \end{equation*}

    \item For any pullback square
      \begin{equation*}
        \begin{tikzcd}
          y'
          \arrow[r, "f'"]
          \arrow[d, swap, "g'"]
          & x'
          \arrow[d, "g"]
          \\
          y
          \arrow[r, "f"]
          & x
        \end{tikzcd},
      \end{equation*}
      if the arrow $f$ belongs to $\category{C}_{\dagger}$ and the arrow $g$ belongs to $\category{C}\updag$, then the arrow $f'$ belongs to $\category{C}_{\dagger}$ and the arrow $g'$ belongs to $\category{C}\updag$.\note{This is a slightly weaker condition than that imposed in \cite{spectralmackeyfunctors1}. The results we will need will hold just as well under this weakened assumption.}
  \end{enumerate}
\end{definition}

We will call a square of the form
\begin{equation*}
  \begin{tikzcd}
    y'
    \arrow[r, rightarrowtail]
    \arrow[d, two heads]
    & x'
    \arrow[d, two heads]
    \\
    y
    \arrow[r, rightarrowtail]
    & x
  \end{tikzcd}
\end{equation*}
\emph{ambigressive.} If such an ambigressive square is also a pullback square, we will call it \emph{ambigressive pullback.}

It is now time to set about building our $\infty$-categories of spans. Fix some triple $\triple{C}$. We will define two models for such a category: a Segal space $\SPAN\triple{C}$, and the quasicategory $\Span\triple{C}$.

Given two $1$-simplices in $\Span\triple{C}$ represented by spans $X \twoheadleftarrow Y \rightarrowtail X'$ and $X' \twoheadleftarrow Y' \rightarrowtail X''$ in $\category{C}$, we need to specify what it means to compose them. We define the composition to be the span $X \twoheadleftarrow Y \times_{X'} Y' \rightarrowtail X''$ given by the diagram below, where the top square is pullback.
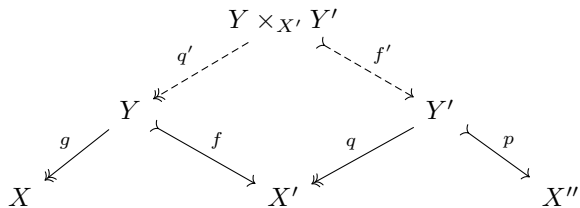
\begin{equation*}
  \begin{tikzcd}
    && Y \times_{X'} Y'
    \arrow[dl, two heads, dashed, swap, "q'"]
    \arrow[dr, tail, dashed, "f'"]
    \\
    & Y
    \arrow[dl, two heads, swap, "g"]
    \arrow[dr, tail, "f"]
    && Y'
    \arrow[dl, two heads, swap, "q"]
    \arrow[dr, tail, "p"]
    \\
    X
    && X'
    && X''
  \end{tikzcd}
\end{equation*}
Note that the conditions in \hyperref[def:adequate_triple]{Definition~\ref*{def:adequate_triple}} tell us precisely that such a composition always exists.

Equivalently, we this means that a $2$-simplex in $\Span\triple{C}$ should be a diagram in $\category{C}$ of the above form. More generally, we would like an $n$-simplex in $\Span\triple{C}$ to correspond to $n$-fold composition of spans. For example, for $n = 4$, such an $n$-simplex in $\Span\triple{C}$ should consist of a diagram in $\category{C}$ of the form
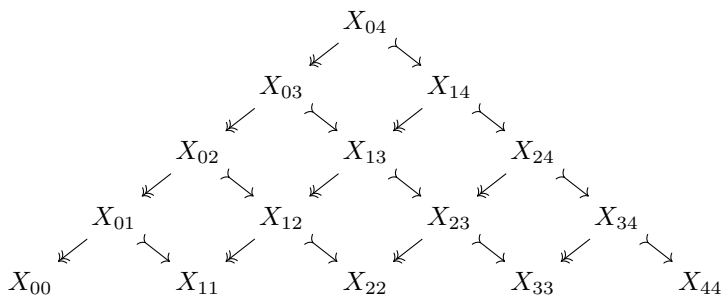
\begin{equation*}
  \begin{tikzcd}[column sep=tiny, row sep=small]
    \\
    &&&& X_{04}
    \arrow[dl, two heads]
    \arrow[dr, tail]
    \\
    &&& X_{03}
    \arrow[dl, two heads]
    \arrow[dr, tail]
    && X_{14}
    \arrow[dl, two heads]
    \arrow[dr, tail]
    \\
    && X_{02}
    \arrow[dl, two heads]
    \arrow[dr, tail]
    && X_{13}
    \arrow[dl, two heads]
    \arrow[dr, tail]
    && X_{24}
    \arrow[dl, two heads]
    \arrow[dr, tail]
    \\
    & X_{01}
    \arrow[dl, two heads]
    \arrow[dr, tail]
    && X_{12}
    \arrow[dl, two heads]
    \arrow[dr, tail]
    && X_{23}
    \arrow[dl, two heads]
    \arrow[dr, tail]
    && X_{34}
    \arrow[dl, two heads]
    \arrow[dr, tail]
    \\
    X_{00}
    && X_{11}
    && X_{22}
    && X_{33}
    && X_{44}
  \end{tikzcd},
\end{equation*}
where the morphisms are ingressive and egressive as shown, and each square is pullback.

All that remains is to make this definition rigorous. To this end, we make the following definition.

\begin{definition}
  for each $n \geq 0$, we define a poset $\Sigma_{n}$ as follows:
  \begin{itemize}
    \item The elements of $\Sigma_{n}$ are pairs of integers $(i, j)$, with $0 \leq i \leq j \leq n$.
  
    \item For $(i, j)$, $(i', j') \in \Sigma_{n}$, we define
      \begin{equation*}
        (i, j) \leq (i', j') \iff i \leq i' \leq j' \leq j.
      \end{equation*}
  \end{itemize}
\end{definition}

\begin{example}
  \label{eg:sigma_4}
  We can draw $\Sigma_{4}$ as follows.
  
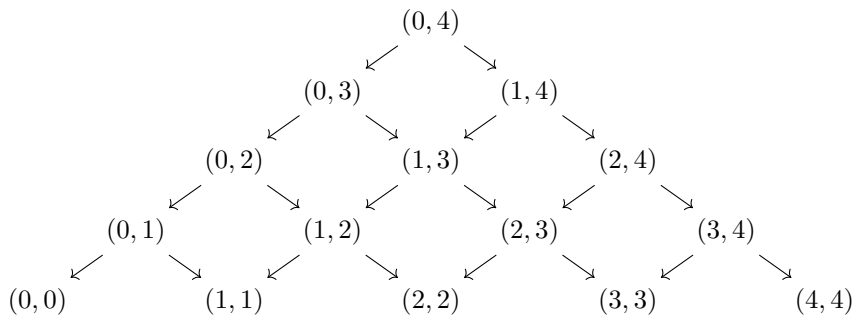
\begin{equation*}
    \begin{tikzcd}[column sep=tiny, row sep=small]
      \\
      &&&& {(0, 4)}
      \arrow[dl]
      \arrow[dr]
      \\
      &&& {(0, 3)}
      \arrow[dl]
      \arrow[dr]
      && {(1, 4)}
      \arrow[dl]
      \arrow[dr]
      \\
      && {(0, 2)}
      \arrow[dl]
      \arrow[dr]
      && {(1, 3)}
      \arrow[dl]
      \arrow[dr]
      && {(2, 4)}
      \arrow[dl]
      \arrow[dr]
      \\
      & {(0, 1)}
      \arrow[dl]
      \arrow[dr]
      && {(1, 2)}
      \arrow[dl]
      \arrow[dr]
      && {(2, 3)}
      \arrow[dl]
      \arrow[dr]
      && {(3, 4)}
      \arrow[dl]
      \arrow[dr]
      \\
      {(0, 0)}
      && {(1, 1)}
      && {(2, 2)}
      && {(3, 3)}
      && {(4, 4)}
    \end{tikzcd}
  \end{equation*}
\end{example}

These posets $\Sigma_{n}$ assemble into a cosimplicial object 
\begin{equation*}
  \Delta \to \SSet;\qquad [n] \mapsto N(\Sigma_{n}),
\end{equation*}
Left Kan extending along the Yoneda embedding $\Delta \hookrightarrow \SSet$ yields a functor $\sd\colon \SSet \to \SSet$. General abstract nonsense gives us right adjoint $\Span'\colon \SSet \to \SSet$ forming an adjunction
\begin{equation*}
  \sd: \SSet \longleftrightarrow \SSet : \Span'.
\end{equation*}
For any simplicial set $A$, the simplicial set $\Span'(A)$ has $n$-simplices
\begin{equation*}
  \Span'(A)_{n} = \Hom_{\SSet}(\sd(\Delta^{n}), A).
\end{equation*}

Note that for a quasicategory $\category{C}$, the simplicial set $\Span'(\category{C})$ is very close to what we want; its $n$-simplices are maps $\Sigma_{n} \to \category{C}$, but the backward- and forward-facing legs of the spans do not necessarily belong to the categories $\category{C}\downdag$ and $\category{C}\updag$, and the squares are not necessarily pullback. We will use the term to refer to a map $\sd(\Delta^{n}) \to \category{C}$ which has the desired form.

%For any $n$-simplex $\sd(\Delta^{n}) \to \category{C}$ corresponding to a diagram as above, we will call the morphisms $X_{ij} \to X_{ij'}$ \emph{backward-facing,} and the morphisms $X_{ij} \to X_{i'j}$ \emph{forward-facing.} 

\begin{definition}
  \label{def:ambigressive_cartesian}
  For any adequate triple $\triple{C}$, we will call a functor $\sd(\Delta^{n}) \to \category{C}$ \emph{ambigressive Cartesian} if each square in $\sd(\Delta^{n}) = N(\Sigma_{n})$ of the form
  \begin{equation*}
    \begin{tikzcd}[column sep=tiny]
      & {(i, j)}
      \arrow[dr]
      \arrow[dl]
      \\
      {(i, j - \ell)}
      \arrow[dr]
      & & {(i + k , j)}
      \arrow[dl]
      \\
      & {(i + k , j - \ell)}
    \end{tikzcd}
  \end{equation*}
  (where we include the possibilities $k = 0$ and $\ell = 0$) is mapped to an ambigressive pullback square,
  \begin{equation*}
    \begin{tikzcd}[column sep=tiny]
      & X_{ij}
      \arrow[dr]
      \arrow[dl]
      \\
      X_{i(j-\ell)}
      \arrow[dr]
      & & X_{(i+k)(j)}
      \arrow[dl]
      \\
      & X_{(i+k)(j-\ell)}
    \end{tikzcd}
  \end{equation*}
  where the backwards-facing morphsims are egressive and the forwards-facing morphisms are ingressive. 
\end{definition}

We are now ready to define our category $\Span(\category{C})$.

\begin{definition}
  We define a simplicial set $\Span\triple{C}$ level-wise to be the subset
  \begin{equation*}
    \Span\triple{C}_{n} \subseteq \Span'(\category{C})_{n}
  \end{equation*}
  on functors $\sd(\Delta^{n}) \to \category{C}$ which are ambigressive Cartesian.
\end{definition}

While it is not hard to check that the face and degeneracy maps on $\Span'\triple{C}$ restrict to $\Span\triple{C}$ (and thus that the above construction really produces a simplicial set), solving the lifting problems necessary to prove directly that $\Span\triple{C}$ is a quasicategory turns out to be combinatorially strenuous. It turns out to be easier to switch to a new model of quasicategories, (complete) Segal spaces, which give easier access to homotopical data.

To that end, for an adequate triple $\triple{C}$, we make the following defintion.
\begin{definition}
  \label{def:acart}
  For any simplicial set $A$, we define a simplicial set $\Map^{\mathrm{aCart}}(\sd(A), \category{C})^{\simeq}$ to be the full simplicial subset
  \begin{equation*}
    \Map^{\mathrm{aCart}}(\sd(A), \category{C})^{\simeq} \subseteq \Map(\sd(A), \category{C})^{\simeq}
  \end{equation*}
  on functors $\sd(A) \to \category{C}$ such that for each standard simplex $\sigma$ of $A$, the image of $\sd(\sigma)$ in $\category{C}$ is ambigressive Cartesian. 
\end{definition}

\begin{definition}
  The (complete) Segal space of spans $\SPAN\triple{C}$ is defined level-wise by 
  \begin{equation*}
    \SPAN\triple{C}_{n} = \Map^{\mathrm{aCart}}(\sd(\Delta^{n}), \category{C})^{\simeq}.
  \end{equation*}
\end{definition}

For a proof that this really is a Segal space, the reader is once again referred to \cite{spectralmackeyfunctors1}. Note that the quasicategory of spans $\Span\triple{C}$ is the bottom row of the complete Segal space of span $\SPAN\triple{C}$; that is,
\begin{equation*}
  \Span\triple{C} = \SPAN\triple{C}/\Delta^{0}.
\end{equation*}
\hyperref[prop:joyal-tierney_quillen_equivalence]{Proposition~\ref*{prop:joyal-tierney_quillen_equivalence}} thus immediately implies that $\Span\triple{C}$ is a quasicategory.

\subsection{Important example: spans and equivalences}
\label{ssc:important_example_spans_and_equivalences}

Let $\category{C}$ be a quasicategory. One can think about spans in $\category{C}$ as `fractions' of morphisms in $\category{C}$. Consider a span in $\category{C}$
\begin{equation*}
  \begin{tikzcd}
    & Y
    \arrow[dl, swap, "g"]
    \arrow[dr, "f"]
    \\
    X
    && X'
  \end{tikzcd}.
\end{equation*}
If the morphism $g$ were invertible, we could fill this to a full 2-simplex
\begin{equation*}
  \begin{tikzcd}
    & Y
    \arrow[dl, swap, "g"]
    \arrow[dr, "f"]
    \\
    X
    \arrow[rr, swap, dashed, "f \circ g^{-1}"]
    && X'
  \end{tikzcd}.
\end{equation*}
Even if $g$ is not invertible, we can think of our span as formally representing a `fraction of morphisms' $f \circ g^{-1}\colon X \to X'$, which may or may not exist in $\category{C}$.

With this point of view in mind, if we build from a quasicategory $\category{C}$ a category of spans whose forwards-facing legs are arbitrary morphisms in $\category{C}$, and whose backward-facing legs are equivalences, then we have not really done anything; each span, thought of as a formal morphism in $\category{C}$, is represented by an actual morphism in $\category{C}$. Thus, the category of spans in $\category{C}$ whose backward-facing legs are equivalences should be equivalent to the category $\category{C}$. In the following example, we construct this equivalence.

\begin{example}
  \label{eg:spans_with_equivalences_on_one_leg}
  For any quasicategory $\category{C}$, we can define a triple $(\category{C}, \category{C}\updag = \category{C}, \category{C}\downdag = \category{C}^{\simeq})$. These choices of ingressive and egressive morphisms correspond to spans of the form
  \begin{equation*}
    \begin{tikzcd}
      & Y
      \arrow[dl, swap, "\simeq"]
      \arrow[dr]
      \\
      X
      && X'
    \end{tikzcd},
  \end{equation*}
  i.e.\ spans such that the backwards-facing map $Y \to X$ is an equivalence in $\category{C}$. With this triple, a functor $\sd(\Delta^{n}) \to \category{C}$ is ambigressive cartesian if and only if each backward-facing leg is an equivalence. Note that the triple $(\category{C}, \category{C}, \category{C}^{\simeq})$ defined above is adequate even if $\category{C}$ does not admit pullbacks: for any solid diagram below, where $g$ is an equivalence, the dashed square completion is a pullback, where $g^{-1}$ is any homotopy inverse for $g$.
  \begin{equation*}
    \begin{tikzcd}
      X
      \arrow[r, dashed, "g^{-1} \circ f"]
      \arrow[d, swap, dashed, "\id"]
      & X'
      \arrow[d, "g"]
      \\
      X
      \arrow[r, "f"]
      & Y'
    \end{tikzcd}
  \end{equation*}

  We will denote the resulting Segal space of spans by $\SPAN^{\simeq}(\category{C})$, and the quasicategory $\SPAN^{\simeq}(\category{C})/\Delta^{0}$ by $\Span^{\simeq}(\category{C})$. 

  The quasicategory $\Span^{\simeq}(\category{C})$ is categorically equivalent to the category $\category{C}$ from which we started. More specifically, there exists a weak categorical equivalence $\category{C} \to \Span^{\simeq}(\category{C})$. We will now construct this map. 

  For each $n \geq 0$ there is a retraction of posets
  \begin{equation*}
    [n] \hookrightarrow \Sigma_{n} \twoheadrightarrow [n],
  \end{equation*}
  where the first map takes $i \mapsto (i, n)$, and the second takes $(i, j) \mapsto i$. In the case $n = 2$, the inclusion can be drawn

  \tikzset{LA/.style = {line width=#1, -{Straight Barb[length=3pt]}}, LA/.default=1.5pt}
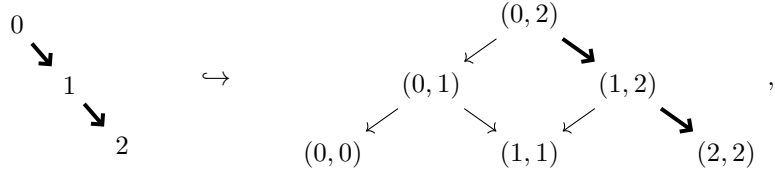
  \begin{equation*}
    \begin{tikzcd}[row sep=small, column sep=tiny]
      0
      \arrow[dr, LA]
      \\
      & 1
      \arrow[dr, LA]
      \\
      && 2
    \end{tikzcd}
    \qquad \hookrightarrow \qquad
    \begin{tikzcd}[row sep=small, column sep=tiny]
      && (0, 2)
      \arrow[dl]
      \arrow[dr, LA]
      \\
      & (0, 1)
      \arrow[dl]
      \arrow[dr]
      && (1, 2)
      \arrow[dl]
      \arrow[dr, LA]
      \\
      (0, 0)
      && (1, 1)
      && (2, 2)
    \end{tikzcd},
  \end{equation*}
  and the map $\sd(\Delta^{n}) \to \Delta^{n}$ `collapses' the backward-facing legs of the spans. The general case is analogous.

  Applying the nerve, we find a retraction of simplicial sets
  \begin{equation*}
    \Delta^{n} \hookrightarrow \sd(\Delta^{n}) \twoheadrightarrow \Delta^{n}.
  \end{equation*}
  Pulling back along these maps gives us, for each $n$, a retraction
  \begin{equation*}
    \begin{tikzcd}
      \Map(\Delta^{n}, \category{C})^{\simeq}
      \arrow[r, hook, "i'_{n}"]
      & \Map(\sd(\Delta^{n}), \category{C})^{\simeq} 
      \arrow[r, hook, "r'_{n}"]
      & \Map(\Delta^{n}, \category{C})^{\simeq}
    \end{tikzcd}.
  \end{equation*}
  Note that the image of the inclusion $i'_{n}$ consists of spans whose backwards-facing legs are identities (hence certainly equivalences), and hence that we can restrict the total space of our retraction:
  \begin{equation*}
    \begin{tikzcd}
      \Map(\Delta^{n}, \category{C})^{\simeq} 
      \arrow[r, hook, "i_{n}"]
      & \Map^{\aCart}(\sd(\Delta^{n}), \category{C})^{\simeq}
      \arrow[r, two heads, "r_{n}"]
      & \Map(\Delta^{n}, \category{C})^{\simeq}
    \end{tikzcd}.
  \end{equation*}

  It is not hard to convince oneself that a map $\sd(\Delta^{n}) \to \category{C}$ is ambigressive cartesian if and only if it is a right Kan extension of its restriction along $\Delta^{n} \hookrightarrow \sd(\Delta^{n})$. In the case $n = 2$, for example, it is clear from the limit formula that any right Kan extension of a diagram $\Delta^{2} \to \category{C}$ along the inclusion $\Delta^{2} \hookrightarrow \sd(\Delta^{2})$ is of the form
  
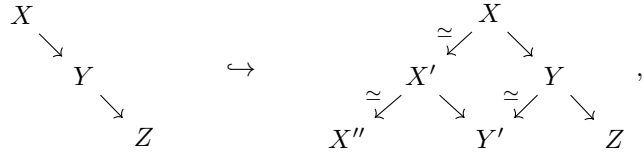
\begin{equation*}
    \begin{tikzcd}[row sep=small, column sep=tiny]
      X
      \arrow[dr]
      \\
      & Y
      \arrow[dr]
      \\
      && Z
    \end{tikzcd}
    \qquad \hookrightarrow \qquad
    \begin{tikzcd}[row sep=small, column sep=tiny]
      && X
      \arrow[dl, swap, "\simeq"]
      \arrow[dr]
      \\
      & X'
      \arrow[dl, swap, "\simeq"]
      \arrow[dr]
      && Y
      \arrow[dl, swap, "\simeq"]
      \arrow[dr]
      \\
      X''
      && Y'
      && Z
    \end{tikzcd},
  \end{equation*}
  (and conversely, any diagram of the above form is a right Kan extension), and diagrams of the above form are precisely those which are ambigressive cartesian. Put differently, $\Fun^{\aCart}(\sd(\Delta^{n}), \category{C})$ is the full subcategory of $\Fun(\sd(\Delta^{n}), \category{C})$ on functors which are right Kan extensions of their restrictions to $\Delta^{n}$. This, together with \cite[Prop.\ 4.3.2.15]{highertopostheory} implies that $r_{n}$ is a trivial fibration for all $n$. The identity $\id_{\Fun(\Delta^{n}, \category{C})}$ is certainly a weak Kan equivalence, so by the $2/3$ property for weak equivalences, $i_{n}$ is also a weak equivalence for all $n$.

  The projections $\sd(\Delta^{n}) \to \Delta^{n}$ assemble into a cosimplicial object $\sd(\Delta^{\bullet}) \to \Delta^{\bullet}$. From this it easily follows that the maps $i_{n}$ assemble into a map of Segal spaces
  \begin{equation*}
    i\colon \Gamma(\category{C}) = \Fun(\Delta^{\bullet}, \category{C}) \to \Fun^{\aCart}(\sd(\Delta^{\bullet}), \category{C}) = \SPAN^{\simeq}(\category{C})
  \end{equation*}
  which is a level-wise Kan weak equivalence, and hence a weak equivalence in the complete Segal spaces model structure.\footnote{For a brief review of the model structures in question, see \hyperref[sss:a_brief_list_of_model_structures]{Appendix~\ref*{sss:a_brief_list_of_model_structures}}. For more information, see \cite{qcats_vs_segal_spaces}, especially the discussion following Thm.\ 4.1.} Here $\Gamma(\category{C})$ is the complete Segal space defined in \cite{qcats_vs_segal_spaces} in the discussion preceding Thm.\ 4.11. Thus, because $\Gamma(\category{C})$ and $\SPAN^{\simeq}(\category{C})$ are complete Segal spaces (hence fibrant-cofibrant objects in the complete Segal space model structure), the restricion of $i$ to first rows $i / \Delta^{0}$ is a weak categorical equivalence, which is what we wanted to show.
\end{example}

As \hyperref[eg:spans_with_equivalences_on_one_leg]{Example~\ref*{eg:spans_with_equivalences_on_one_leg}} shows, taking spans whose backward-facing legs are equivalences allows us to find an equivalence between $\category{C}$ and a quasicategory of spans in $\category{C}$. Similarly, if we look at spans whose forward-facing legs are equivalences, we find a category of spans equivalent to $\category{C}\op$.

\begin{example}
  \label{eg:spans_with_equivalences_on_forwards_leg}
  For any quasicategory $\category{C}$, we can define an adequate triple $(\category{C}, \category{C}\downdag = \category{C}^{\simeq}, \category{C}\updag = \category{C})$. These selections of ingressive and egressive morphisms correspond to spans of the form
  \begin{equation*}
    \begin{tikzcd}
      & Y
      \arrow[dl]
      \arrow[dr, "\simeq"]
      \\
      X
      && X'
    \end{tikzcd},
  \end{equation*}
  i.e.\ spans such that the forward-facing map $Y \to X$ is an equivalence in $\category{C}$. We will denote the corresponding category of spans by $\Span_{\simeq}(\category{C})$. Exactly analogous reasoning to that in \hyperref[eg:spans_with_equivalences_on_one_leg]{Example~\ref*{eg:spans_with_equivalences_on_one_leg}} gives us a categorical equivalence
  \begin{equation*}
    \category{C}\op \to \Span_{\simeq}(\category{C}).
  \end{equation*}
\end{example}

\subsection{Maps between triples}
\label{ssc:maps_between_triples}

We have now seen that given any quasicategory $\category{C}$ with pullbacks, one can create a quasicategory of spans $\Span(\category{C})$, and that by defining an adequate triple structure $\triple{C}$ on $\category{C}$, one can specify certain classes of morphisms to which the legs of the spans are allowed to belong. It is easy to see that any functor $\category{C} \to \category{D}$ between quasicategories $\category{C}$ and $\category{D}$ with pullbacks induces a functor $\Span(\category{C}) \to \Span(\category{D})$. Similarly, we make the following definition.

\begin{definition}
  Let $\triple{C}$ and $\triple{D}$ be adequate triples. A functor $p\colon \category{C} \to \category{D}$ is said to be a \defn{functor between adequate triples} $\triple{C} \to \triple{D}$ 
  if the following conditions hold. 
  \begin{enumerate}
    \item The functor $p$ preserves ingressive morphisms. That is, $p(\category{C}_{\dagger}) \subseteq \category{D}_{\dagger}$

    \item The functor $p$ preserves egressive morphsims. That is, $p(\category{C}^{\dagger}) \subseteq \category{D}^{\dagger}$. 

    \item The functor $p$ preserves ambigressive pullbacks.
  \end{enumerate}
\end{definition}

We will denote the restriction of $p$ to $\category{C}_{\dagger}$ by
\begin{equation*}
  p_{\dagger}\colon\category{C}_{\dagger} \to \category{D}_{\dagger},
\end{equation*}
and similarly for $p^{\dagger}$. Clearly, a functor $p$ between adequate triples gives functors $\Span(p)$ and $\SPAN(p)$ respectively between quasicategories and Segal spaces of spans.

It is natural to wonder about the relationship between $p$ and $\Span(p)$. The goal of the remainder of this section is to provide a proof of the following theorem, originally appearing as \cite[Thm.\ 12.2]{spectralmackeyfunctors1}. This theorem provides conditions on a functor $p$ between adequate triples such that the induced map between category of spans is an inner fibration, and provides the form of a $p$-cocartesian morphism in $\Span\triple{C}$.

\begin{theorem}
  \label{thm:main}
  Let $p\colon \triple{C} \to \triple{D}$ be a functor between adequate triples such that $p\colon \category{C} \to \category{D}$ is an inner fibration which satisfies the following conditions.
  \begin{enumerate}
    \item Each morphism $g \in \category{D}\downdag$ admits a lift to a morphism in $\category{C}\downdag$ (given a lift of the source) which is both $p$-cocartesian and $p\downdag$-cocartesian.

    \item Consider a commutative square
      \begin{equation*}
        \sigma = \quad
        \begin{tikzcd}
          y'
          \arrow[r, rightarrowtail, "f'"]
          \arrow[d, two heads, swap, "g'"]
          & x'
          \arrow[d, "g"]
          \\
          y
          \arrow[r, rightarrowtail, "f"]
          & x
        \end{tikzcd}
      \end{equation*}
      in $\category{C}$ such that $p(\sigma)$ is ambigressive pullback in $\category{D}$, where $g'$ belongs to $\category{C}\updag$, and $f$ and $f'$ belong to $\category{C}\downdag$. Suppose that $f$ is $p$-cocartesian. Then $f'$ is $p'$-cocartesian if and only if $\sigma$ is an ambigressive pullback square (and in particular $g \in \category{C}\updag$).
  \end{enumerate}
  Then $\pi$ is an inner fibration. Further, if a morphism in $\Span\triple{C}$ is of the form
  \begin{equation*}
    \begin{tikzcd}
      & y
      \arrow[dl, two heads, swap, "\phi"]
      \arrow[dr, rightarrowtail, "\psi"]
      \\
      x
      && x'
    \end{tikzcd},
  \end{equation*}
  where $\phi$ is egressive and $p^{\dagger}$-cartesian, and $\psi$ is ingressive and $p$-cocartesian, it is $\pi$-cocartesian.%\footnote{This is not precisely the same condition that Barwick gives: he demands that $\phi$ be $p$-cartesian rather than $p^{\dagger}$-cartesian. As far as I can tell, this is simply an error; Barwick seems to use the condition as stated above in his proof.}
\end{theorem}

This theorem originally appears in \cite{spectralmackeyfunctors1}. However, the proof there is combinatorial in nature. Our proof uses the material developed in \hyperref[sec:cocartesian_fibrations_between_complete_segal_spaces]{Section~\ref*{sec:cocartesian_fibrations_between_complete_segal_spaces}}, leveraging the fact that $\infty$-categories of spans have natural incarnations as (complete) Segal spaces.

Before we get to the meat of the proof, however, we have to prove some preliminary results.

We also note a variant, the proof of which is similar but easier.

\begin{theorem}
  \label{thm:new_barwick}
  Let $p\colon \triple{C} \to \triple{D}$ be a functor between adequate triples such that $p\colon \category{C} \to \category{D}$ is an inner fibration which satisfies the following conditions.
  \begin{enumerate}
    \item The subcategory $\category{C}\updag \subseteq \category{C}$ consists of all $p$-cartesian morphisms in $\category{C}$; that is, an $n$-simplex in $\category{C}$ belongs to $\category{C}\updag$ if and only if each $1$-simplex it contains is $p$-cartesian.

    \item The map $p\updag\colon \category{C}\updag \to \category{D}\updag$ is a cartesian fibration.

    \item Consider a square
      \begin{equation*}
          \sigma = \quad
          \begin{tikzcd}
            y'
            \arrow[r, "f'"]
            \arrow[d, two heads, swap, "g'"]
            & x'
            \arrow[d, two heads, "g"]
            \\
            y
            \arrow[r, rightarrowtail, "f"]
            & x
          \end{tikzcd}
      \end{equation*}
      in $\category{C}$ where $g$ and $g'$ belong to $\category{C}\updag$, and $f$ belongs to $\category{C}\downdag$. Further suppose that $f$ is $p$-cocartesian. Then $f'$ belongs to $\category{C}\downdag$, and is both $p$-cocartesian and $p\downdag$-cocartesian.
  \end{enumerate}
  Then spans of the form
  \begin{equation*}
    \begin{tikzcd}
      & z
      \arrow[dl, two heads, swap, "g"]
      \arrow[dr, rightarrowtail, "f"]
      \\
      x
      && y
    \end{tikzcd}
  \end{equation*}
  are cocartesian, where $g$ is $p\updag$-cartesian and $f$ is $p$-cocartesian.
\end{theorem}
%\begin{proof}
%  We first note that $\sigma$ is automatically pullback; since the sides are $p$-cartesian, it is automatically a relative pullback, and it lies over a pullback square.
%
%  The rest is even simpler than the proof in my thesis of the ordinary Barwick's theorem. The square that we have to show is homotopy pullback factors as before, but this time we don't have to take the final homotopy pullback to find path components corresponding to cartesian 2-simplices, since our 2-simplices are automatically cartesian.
%\end{proof}

\subsection{Housekeeping}

In this section we recall some basic facts we will need in the course of our proof of the main theorem. It is recommended to skip this section on first reading, and refer back to the results as necessary.

\subsubsection{Isofibrations}

We will need some results about isofibrations, sometimes called categorical fibrations. These are standard, and proofs can be found in \cite{kerodon}.

\begin{proposition}
  \label{prop:characterization_of_isofibrations}
  Let $f\colon \category{C} \to \category{D}$ be a map between quasicategories. The following are equivalent:
  \begin{itemize}
    \item The map $f$ is an isofibration.

    \item The map $f$ is an inner fibration, and each equivalence in $\category{D}$ has a lift in $\category{C}$ which is an equivalence.
  \end{itemize}
\end{proposition}

We will denote by $\Map(-, -)$ the internal hom on $\SSet$.
\begin{proposition}
  Let $f\colon \category{C} \to \category{D}$ be an isofibration between quasicategories, and let $i\colon A \hookrightarrow B$ be an inclusion of simplicial sets. Then the map
  \begin{equation*}
    \Map(B, \category{C}) \to \Map(A, \category{C}) \times_{\Map(A, \category{D})} \Map(B, \category{D})
  \end{equation*}
  is an isofibration.
\end{proposition}

We will denote by $(-)^{\simeq}$ the \emph{core} functor, i.e.\ the functor which takes a simplicial set $X$ to the largest Kan complex it contains.
\begin{proposition}
  Let $f\colon X \to Y$ be an isofibration of simplicial sets. Then
  \begin{equation*}
    f^{\simeq}\colon X^{\simeq} \to Y^{\simeq}
  \end{equation*}
  is a Kan fibration.
\end{proposition}

We will denote the core of the functor $\Map(-, -)$ by $\Fun(-, -)$. That is, for simplicial sets $A$ and $B$, we have
\begin{equation*}
  \Fun(A, B) \cong \Map(A, B)^{\simeq}.
\end{equation*}

\begin{corollary}
  \label{cor:kan_fib_and_isofib}
  Let $f\colon \category{C} \to \category{D}$ be an isofibration between quasicategories, and let $A \hookrightarrow B$ be an inclusion of simplicial sets. Then the map
  \begin{equation*}
    \Fun(B, \category{C}) \to \Fun(A, \category{C}) \times_{\Fun(A, \category{D})} \Fun(B, \category{D})
  \end{equation*}
  is a Kan fibration.
\end{corollary}

\subsubsection{Connected components}

Our proof will rely heavily on the fact that the many properties that morphisms can have (ingressive, egressive, $p$-cocartesian, etc.) are well-behaved with respect to the homotopical structure of the complete Segal spaces in which they live; for example, if a morphism is $p$-cocartesian, then every morphism in its path component is $p$-cocartesian.

Let $X'$ and $X$ be simplicial sets, and let $f\colon X' \hookrightarrow X$ be an inclusion. We will say that $f$ is an \emph{inclusion of connected components} if for all simplices $\sigma \in X$, if $\sigma$ has any vertex in common with $X'$, then $\sigma$ is wholly contained in $X'$.

\begin{lemma}
  \label{lemma:lifting_wrt_inclusions_of_connected_components}
  Let $A$ be a simplicial set such that between any two vertices $x$ and $y$ of $A$ there exists a finite zig-zag of 1-simplices of $A$ connecting $x$ to $y$. Let $A_{0} \subseteq$ A be a nonempty simplicial subset of $A$. Let $f\colon X' \hookrightarrow X$ be a morphism between simplicial sets which is an inclusion of connected components. Then the following dashed lift always exists.
  \begin{equation*}
    \begin{tikzcd}
      A_{0}
      \arrow[r, "f_{0}"]
      \arrow[d, hook]
      & X'
      \arrow[d, hook]
      \\
      A
      \arrow[r, "f"]
      \arrow[ur, dashed]
      & X
    \end{tikzcd}
  \end{equation*}
\end{lemma}
\begin{proof}
  We have a map $f\colon A \to X$; in order to construct our dashed lift, it suffices to show that under the above assumptions, $f$ takes every simplex of $A$ to a simplex which belongs to $X'$. To this end, let $\sigma \in A$ be a simplex. There exists a finite zig-zag of $1$-simplices connecting some vertex of $\sigma$ to a vertex $\alpha$ of $A_{0}$; under $f$, this is mapped to a zig-zag of $1$-simplices connecting $f(\alpha)$ to $f(\sigma)$. The vertex $f(\alpha)$ belongs to $X'$; proceeding inductively, we find that the image of every 1-simplex belonging to the zig-zag also belongs to $X'$, and that $f(\sigma)$ therefore also belongs to $X'$.
\end{proof}

\begin{corollary}
  \label{cor:connected_components_kan_fibration}
  Given any commuting square of simplicial sets
  \begin{equation*}
    \begin{tikzcd}
      X'
      \arrow[r, hook, "i"]
      \arrow[d, swap, "f'"]
      & X
      \arrow[d, "f"]
      \\
      Y'
      \arrow[r, hook]
      & Y
    \end{tikzcd}
  \end{equation*}
  with monomorphisms as marked, where $i$ is an inclusion of connected components between Kan complexes and $f$ is a Kan fibration, the map $f'$ is a Kan fibration.

  Furthermore, if $f$ is a trivial fibration and $f'$ is surjective on vertices, then $f'$ is a trivial fibration.
\end{corollary}
\begin{proof}
  First, we show that if $f$ is a Kan fibration, then $f'$ is a Kan fibration. We need to show that a dashed lift below exists for each $n \geq 1$ and each $0 \leq i \leq n$.
  \begin{equation*}
    \begin{tikzcd}
      \Lambda^{n}_{i}
      \arrow[r]
      \arrow[d]
      & X'
      \arrow[r, hook, "i"]
      \arrow[d, "f'"]
      & X
      \arrow[d, "f"]
      \\
      \Delta^{n}
      \arrow[r]
      \arrow[ur, dashed]
      & Y'
      \arrow[r, hook]
      & Y
    \end{tikzcd}
  \end{equation*}
  We can always solve our outer lifting problem, and by \hyperref[lemma:lifting_wrt_inclusions_of_connected_components]{Lemma~\ref*{lemma:lifting_wrt_inclusions_of_connected_components}}, our lift factors through $X'$.

  Now suppose that $f$ is a trivial Kan fibration. The logic above says that $f'$ has the right lifting property with respect to $\partial \Delta^{n} \hookrightarrow \Delta^{n}$ for $n \geq 1$; the right lifting property with respect to $\partial \Delta^{0} \hookrightarrow \Delta^{0}$ is equivalent to the surjectivity of $f'$ on vertices.
\end{proof}

\subsubsection{Contractibility}

We will need at certain points to show that various spaces of lifts are contractible. These are mostly common-sense results, but it will be helpful to have them written down somewhere so we can refer to them later.

\begin{lemma}
  \label{lemma:inner_fib_contractible_fibers}
  Let $f\colon \category{C} \to \category{D}$ be an inner fibration between quasicategories. Then for any commuting square
  \begin{equation*}
    \begin{tikzcd}
      \Lambda^{2}_{1}
      \arrow[r, "\alpha"]
      \arrow[d, hook]
      & \category{C}
      \arrow[d, "f"]
      \\
      \Delta^{2}
      \arrow[r, "\beta"]
      & \category{D}
    \end{tikzcd},
  \end{equation*}
  the fiber $F$ in the pullback square below is contractible.
  \begin{equation*}
    \begin{tikzcd}[column sep=large]
      F
      \arrow[r]
      \arrow[d]
      & \Fun(\Delta^{2}, \category{C})
      \arrow[d]
      \\
      \Delta^{0}
      \arrow[r, "{(\alpha, f\alpha, \beta)}"]
      & \Fun(\Lambda^{2}_{1}, \category{C}) \times_{\Fun(\Lambda^{2}_{1}, \category{D})} \Fun(\Delta^{2}, \category{D})
    \end{tikzcd}
  \end{equation*}
\end{lemma}
\begin{proof}
  The right-hand map is a trivial Kan fibration.
\end{proof}

One should interpret this as telling us that given a $\Lambda^{2}_{2}$-horn $\alpha$ in $\category{C}$ lying over a 2-simplex $\beta$ in $\category{D}$, the space of fillings of $\alpha$ lying over $\beta$ is contractible. We will need a pantheon of similar results.

\begin{lemma}
  \label{lemma:cartesian_horn_fillings_contractible}
  Let $f\colon \category{C} \to \category{D}$ be an inner fibration between quasicategories. Let $e$ be any $f$-cartesian edge. Then for any square
  \begin{equation*}
    \begin{tikzcd}
      \Lambda^{2}_{2}
      \arrow[r, "\alpha"]
      \arrow[d, hook]
      & \category{C}
      \arrow[d, "f"]
      \\
      \Delta^{2}
      \arrow[r, "\beta"]
      & \category{D}
    \end{tikzcd}
  \end{equation*}
  such that $\alpha|_{\{1,2\}} = e$, the fiber $F$ in the pullback diagram
  \begin{equation}
    \label{eq:fiber_of_inn_fib_with_cocart}
    \begin{tikzcd}[column sep=large]
      F
      \arrow[r]
      \arrow[d]
      & \Fun(\Delta^{2}, \category{C})
      \arrow[d]
      \\
      \Delta^{0}
      \arrow[r, "{(\alpha, f\alpha, \beta)}"]
      & \Fun(\Lambda^{2}_{2}, \category{C}) \times_{\Fun(\Lambda^{2}_{1}, \category{D})} \Fun(\Delta^{2}, \category{D})
    \end{tikzcd}
  \end{equation}
  is contractible.
\end{lemma}
\begin{proof}
  Denote $f(e) = \bar{e}$. Define a marking $\mathcal{E}'$ on $\category{D}$ containing all degenerate edges and $\bar{e}$. Define a marking $\mathcal{E}$ on $\category{C}$ containing the $f$-cocartesian lifts of the edges in $\mathcal{E}$. By \cite[Prop.\ 3.1.1.6]{highertopostheory}, the map $\category{C}^{\mathcal{E}} \to \category{D}^{\mathcal{E}'}$ has the right-lifting property with respect to all cocartesian-marked anodyne morphisms. It follows that
  \begin{equation*}
    \begin{tikzcd}[column sep=large]
      \Map^{\sharp}((\Delta^{2})^{\mathcal{L}}, \category{C}^{\mathcal{E}})
      \arrow[d]
      \\
      \Map^{\sharp}((\Lambda^{2}_{2})^{\mathcal{L}}, \category{C}^{\mathcal{E}}) \times_{\Map^{\sharp}((\Lambda^{2}_{1})^{\mathcal{L}}, \category{D}^{\mathcal{E}'})} \Map^{\sharp}((\Delta^{2})^{\mathcal{L}}, \category{D}^{\mathcal{E}'})
    \end{tikzcd}
  \end{equation*}
  is a trivial Kan fibration; the map on the right-hand side of \hyperref[eq:fiber_of_inn_fib_with_cocart]{Diagram~\ref*{eq:fiber_of_inn_fib_with_cocart}} is the core of this map, and is thus also a trivial Kan fibration.
\end{proof}

\begin{lemma}
  \label{lemma:twice_cartesian_lifts_contractible}
  Let $\triple{C} \to \triple{D}$ be an inner fibration between triples satisfying the conditions of \hyperref[thm:main]{Theorem~\ref*{thm:main}}. Denote by $\Fun'(\Delta^{1}, \category{C}^{\dagger})$ the full simplicial subset of $\Fun(\Delta^{1}, \category{C}^{\dagger})$ on edges $\Delta^{1} \to \category{C}^{\dagger}$ which are both $p$-cocartesian and $p^{\dagger}$-cocartesian. Then the fibers of the map
  \begin{equation*}
    \Fun'(\Delta^{1}, \category{C}^{\dagger}) \to
    \Fun(\Delta^{\{0\}}, \category{C}^{\dagger})
    \times_{\Fun(\Delta^{\{0\}}, \category{D}^{\dagger})}
    \Fun(\Delta^{1}, \category{D}^{\dagger})
  \end{equation*}
  are contractible.
\end{lemma}
\begin{proof}
  Denote by $\Fun''(\Delta^{1}, \category{C}^{\dagger})$ the full simplicial subset of $\Fun(\Delta^{1}, \category{C}^{\dagger})$ on edges which are $p^{\dagger}$-cartesian. Then we have a square
  \begin{equation*}
    \noindent\makebox[\textwidth]{%
      \begin{tikzcd}[ampersand replacement=\&]
        \Fun'(\Delta^{1}, \category{C}^{\dagger})
        \arrow[r, hook]
        \arrow[d]
        \& \Fun''(\Delta^{1}, \category{C}^{\dagger})
        \arrow[d]
        \\
        \Fun(\Delta^{\{0\}}, \category{C}^{\dagger})
        \times_{\Fun(\Delta^{\{0\}}, \category{D}^{\dagger})}
        \Fun(\Delta^{1}, \category{D}^{\dagger})
        \arrow[r, equals]
        \& \Fun(\Delta^{\{0\}}, \category{C}^{\dagger})
        \times_{\Fun(\Delta^{\{0\}}, \category{D}^{\dagger})}
        \Fun(\Delta^{1}, \category{D}^{\dagger})
      \end{tikzcd}
    }
  \end{equation*}
  in which the top map is an inclusion of connected components, the right map is a trivial Kan fibration, and the left map is surjective on vertices. Thus, the left map is a trivial Kan fibration by \hyperref[cor:connected_components_kan_fibration]{Corollary~\ref*{cor:connected_components_kan_fibration}}
\end{proof}

\subsection{Proof of main theorem}

Our goal is now to prove \hyperref[thm:main]{Theorem~\ref*{thm:main}}. As in the statement of the theorem, we fix a functor of triples $p\colon \triple{C} \to \triple{D}$ giving us a functor of quasicategories
\begin{equation*}
  \pi\colon \Span\triple{C} \to \Span\triple{D}.
\end{equation*}
We wish to show that the following hold.
\begin{enumerate}
  \item The map $\pi$ is an inner fibration.

  \item Morphsims of the form
    \begin{equation}
      \label{eq:form_of_p_cocartesian_morphisms}
      \begin{tikzcd}
        & y
        \arrow[dl, two heads, swap, "\phi"]
        \arrow[dr, rightarrowtail, "\psi"]
        \\
        x
        && x'
      \end{tikzcd},
    \end{equation}
    are $\pi$-cocartesian, where $\phi$ is egressive and $p_{\dagger}$-cartesian, and $\psi$ is ingressive and $p$-cocartesian.
\end{enumerate}

The way forward is clear: the functor $\pi$ is the zeroth row of a functor of Segal spaces
\begin{equation*}
  \pi'\colon \SPAN\triple{C} \to \SPAN\triple{D}.
\end{equation*}

The results of \hyperref[sec:cocartesian_fibrations_between_complete_segal_spaces]{Section~\ref*{sec:cocartesian_fibrations_between_complete_segal_spaces}} therefore guarantee the following.
\begin{enumerate}
  \item By \hyperref[cor:reedy_implies_inner]{Corollary~\ref*{cor:reedy_implies_inner}}, in order to show that $\pi$ is an inner fibration, it suffices to show that $\pi'$ is a Reedy fibration.

  \item In order to show that morphisms of the required form are $\pi$-cocartesian, it will suffice to show by \hyperref[cor:cocart_fib_between_css_gives_cocart_fib_of_quasicats]{Corollary~\ref*{cor:cocart_fib_between_css_gives_cocart_fib_of_quasicats}} that morphisms of this form are $\pi'$-cocartesian. 
\end{enumerate}

\subsubsection{The map \texorpdfstring{$\pi$}{pi} is an inner fibration}

We begin by checking that the map $\pi$ is an inner fibration. As stated above, we do this by checking that $\pi'$ is a Reedy fibration. This follows from the following useful fact, implied by the first assumption of \hyperref[thm:main]{Theorem~\ref*{thm:main}}.
\begin{lemma}
  Let $p\colon \triple{C} \to \triple{D}$ be a functor of triples satisfying the conditions of \hyperref[thm:main]{Theorem~\ref*{thm:main}}. Then the map $p\colon \category{C} \to \category{D}$ is an isofibration.
\end{lemma}
\begin{proof}
  By assumption, $p$ is an inner fibration and $\category{C}$ and $\category{D}$ are quasicategories, so it suffices by \hyperref[prop:characterization_of_isofibrations]{Proposition~\ref*{prop:characterization_of_isofibrations}} to show that $p$ admits lifts of equivalences. Each equivalence in $\category{D}$ is in particular ingressive, and hence admits a $p$-cocartesian lift by the assumptions of \hyperref[thm:main]{Theorem~\ref*{thm:main}}; these are automatically equivalences by \cite[Prop.~2.4.1.5]{highertopostheory}.
\end{proof}

\begin{proposition}
  For any functor of adequate triples $p\colon \triple{C} \to \triple{D}$, the map
  \begin{equation*}
    \pi'\colon \SPAN\triple{C} \to \SPAN\triple{D}
  \end{equation*}
  is a Reedy fibration.
\end{proposition}
\begin{proof}
  Consider the following diagram.
  \begin{equation*}
    \noindent\makebox[\textwidth]{%
      \begin{tikzcd}[ampersand replacement=\&]
        \Fun^{\aCart}(\sd(\Delta^{n}), \category{C})
        \arrow[r, hook]
        \arrow[d]
        \& \Fun(\sd(\Delta^{n}), \category{C})
        \arrow[d]
        \\
        \Fun^{\aCart}(\sd(\partial\Delta^{n}), \category{C})
        \times_{\Fun^{\aCart}(\sd(\partial\Delta^{n}), \category{D})}
        \Fun^{\aCart}(\sd(\Delta^{n}), \category{D})
        \arrow[r, hook]
        \& \Fun(\sd(\partial\Delta^{n}), \category{C})
        \times_{\Fun(\sd(\partial\Delta^{n}), \category{D})}
        \Fun(\sd(\Delta^{n}), \category{D})
      \end{tikzcd}.
    }
  \end{equation*}
  The right-hand map is a Kan fibration by \hyperref[cor:kan_fib_and_isofib]{Corollary~\ref*{cor:kan_fib_and_isofib}} (because $\sd$ preserves monomorphisms), and the top map is an inclusion of connected components, so \hyperref[cor:connected_components_kan_fibration]{Corollary~\ref*{cor:connected_components_kan_fibration}} implies that the left-hand map is a Kan fibration, which is what we needed to show.
\end{proof}

\begin{corollary}
  The map $\pi$ is an inner fibration
\end{corollary}
\begin{proof}
  \hyperref[cor:reedy_implies_inner]{Corollary~\ref*{cor:reedy_implies_inner}}.
\end{proof}

\subsubsection{Cocartesian morphisms have the promised form}
\label{sss:cocartesian_morphisms_have_the_promised_form}

\begin{figure}[p]
  \begin{equation*}
    \begin{tikzcd}
      && 11
      \\
      & 01
      \arrow[ur, rightarrowtail, "\bullet" marking]
      \arrow[dl, two heads, "\vartriangleleft" marking]
      \\
      00
      && 02
      \arrow[ll, two heads]
      \arrow[rr, rightarrowtail]
      && 22
    \end{tikzcd}
    \overset{i_{1}}{\longrightarrow}
    \begin{tikzcd}
      && 11
      \\
      & 01
      \arrow[ur, rightarrowtail, "\bullet" marking]
      \arrow[dl, two heads, "\vartriangleleft" marking]
      \\
      00
      && 02
      \arrow[ll, two heads]
      \arrow[rr, rightarrowtail]
      \arrow[ul, two heads]
      && 22
    \end{tikzcd}
  \end{equation*}
  \vspace{2cm}
  \begin{equation*}
    \overset{i_{2}}{\longrightarrow}
    \begin{tikzcd}
      && 11
      \\
      & 01
      \arrow[ur, rightarrowtail, "\bullet" marking]
      \arrow[dl, two heads, "\vartriangleleft" marking]
      \\
      00
      && 02
      \arrow[ll, two heads]
      \arrow[rr, rightarrowtail]
      \arrow[uu]
      \arrow[ul, two heads]
      && 22
    \end{tikzcd}
    \overset{i_{3}}{\longrightarrow}
    \begin{tikzcd}
      && 11
      \\
      & 01
      \arrow[ur, rightarrowtail, "\bullet" marking]
      \arrow[dl, two heads, "\vartriangleleft" marking]
      && 12
      \\
      00
      && 02
      \arrow[ur, sloped, "\bullet \blacktriangleright" marking, rightarrowtail]
      \arrow[ll, two heads]
      \arrow[rr, rightarrowtail]
      \arrow[uu]
      \arrow[ul, two heads]
      && 22
    \end{tikzcd}
  \end{equation*}
  \vspace{2cm}
  \begin{equation*}
    \overset{i_{4}}{\longrightarrow}
    \begin{tikzcd}
      && 11
      \\
      & 01
      \arrow[ur, rightarrowtail, "\bullet" marking]
      \arrow[dl, two heads, "\vartriangleleft" marking]
      && 12
      \arrow[ul]
      \\
      00
      && 02
      \arrow[ur, sloped, "\bullet \blacktriangleright" marking, rightarrowtail]
      \arrow[ll, two heads]
      \arrow[rr, rightarrowtail]
      \arrow[uu]
      \arrow[ul, two heads]
      && 22
    \end{tikzcd}
    \overset{i_{5}}{\longrightarrow}
    \begin{tikzcd}
      && 11
      \\
      & 01
      \arrow[ur, rightarrowtail, "\bullet" marking]
      \arrow[dl, two heads, "\vartriangleleft" marking]
      && 12
      \arrow[ul]
      \arrow[dr, rightarrowtail]
      \\
      00
      && 02
      \arrow[ur, sloped, "\bullet \blacktriangleright" marking, rightarrowtail]
      \arrow[ll, two heads]
      \arrow[rr, rightarrowtail]
      \arrow[uu]
      \arrow[ul, two heads]
      && 22
    \end{tikzcd}
  \end{equation*}
  \caption{A factorization of the inclusion $\sd(\Lambda^{2}_{0}) \hookrightarrow \sd(\Delta^{2})$, where certain morphisms have been labelled.}
  \label{fig:factorization}
\end{figure}
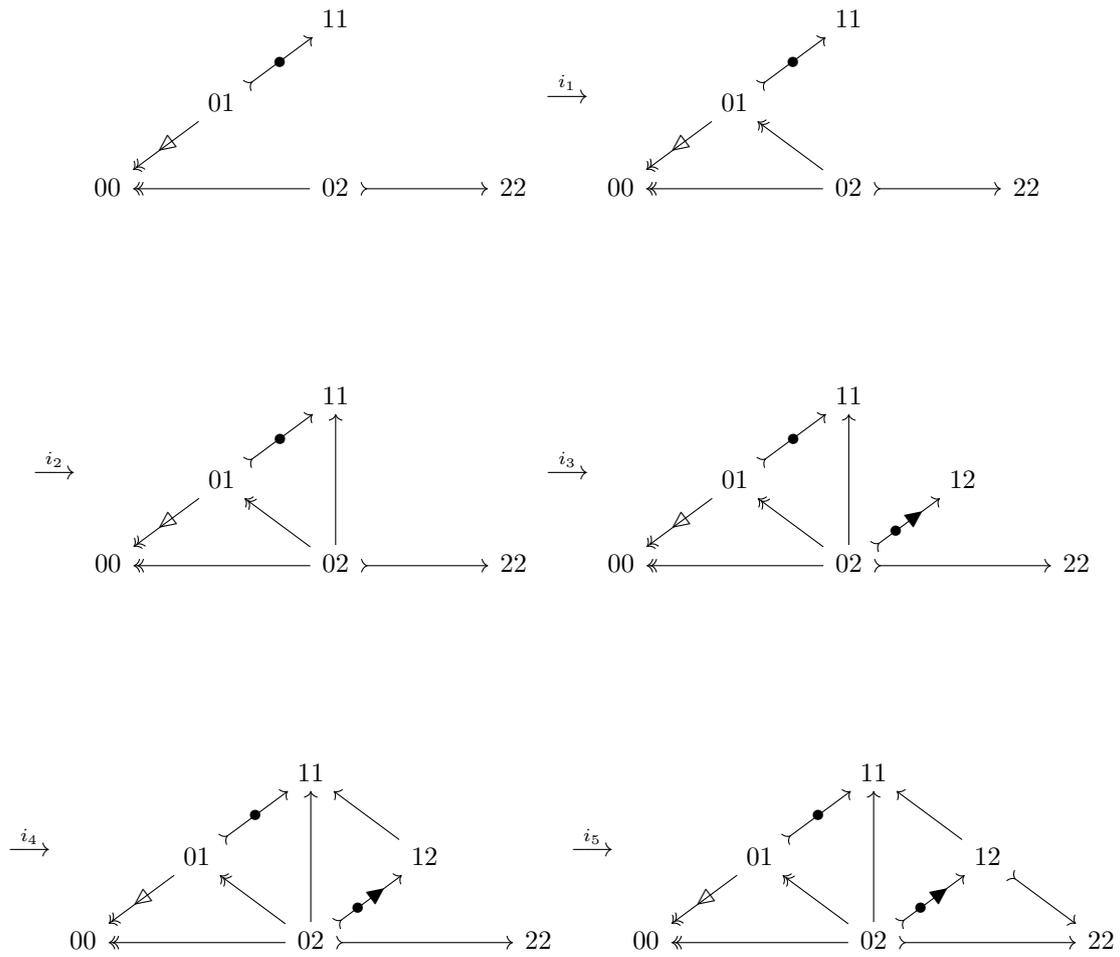

We are now ready to begin in earnest our proof that morphisms of the form \hyperref[eq:form_of_p_cocartesian_morphisms]{Equation~\ref*{eq:form_of_p_cocartesian_morphisms}} are $\pi'$-cocartesian. Morally, this result should not be surprising. We should think of the homotopy pullback condition defining cocartesian morphisms (\hyperref[def:cocartesian_morphism]{Definition~\ref*{def:cocartesian_morphism}}) as telling us that we can fill relative $\Lambda^{2}_{2}$-horns in $\pi'\colon \SPAN\triple{C} \to \SPAN\triple{D}$, and that such fillings are unique up to contractible choice. Finding a filling 
\begin{equation*}
  \begin{tikzcd}
    \Lambda^{2}_{2}
    \arrow[r]
    \arrow[d, hook]
    & \SPAN\triple{C}
    \arrow[d, "\pi'"]
    \\
    \Delta^{2}
    \arrow[r]
    \arrow[ur, dashed]
    & \SPAN\triple{D}
  \end{tikzcd}
\end{equation*}
is equivalent finding a filling
\begin{equation*}
  \begin{tikzcd}
    \sd(\Lambda^{2}_{2})
    \arrow[r]
    \arrow[d, hook]
    & \category{C}
    \arrow[d, "p"]
    \\
    \sd(\Delta^{2})
    \arrow[r]
    \arrow[ur, dashed]
    & \category{D}
  \end{tikzcd}
\end{equation*}
(such that the filling is ambigressive and the necessary square is pullback); the conditions on $p$ guarantee us that we can perform this filling in a series of steps pictured in \hyperref[fig:factorization]{Figure~\ref*{fig:factorization}}, each of which is unique up to contractible choice.

In fact, our proof mainly consists of making this idea rigorous. To this end, we first introduce some notation. Let $p\colon \triple{C} \to \triple{D}$ be a functor of triples.
\begin{itemize}
  \item We will denote $p$-cartesian morphisms in $\category{C}$ with a circle:
    \begin{equation*}
      \begin{tikzcd}
        x
        \arrow[r, "\circ" marking]
        & y
      \end{tikzcd}
    \end{equation*}

  \item We will denote $p$-cocartesian morphisms in $\category{C}$ with a bullet:
    \begin{equation*}
      \begin{tikzcd}
        x
        \arrow[r, "\bullet" marking]
        & y
      \end{tikzcd}
    \end{equation*}

  \item We will denote $p^{\dagger}$-cartesian morphisms in $\category{C}$ with a triangle:
    \begin{equation*}
      \begin{tikzcd}
        x
        \arrow[r, "\vartriangleright" marking]
        & y
      \end{tikzcd}
    \end{equation*}

  \item We will denote $p_{\dagger}$-cocartesian morphisms in $\category{C}$ with a filled triangle:
    \begin{equation*}
      \begin{tikzcd}
        x
        \arrow[r, "\blacktriangleright" marking]
        & y
      \end{tikzcd}
    \end{equation*}
\end{itemize}
Thus, an ingressive morphism which is both $p$-cocartesian and $p_{\dagger}$-cocartesian will be denoted by
\begin{equation*}
  \begin{tikzcd}[column sep=large]
    x
    \arrow[r, rightarrowtail, "\blacktriangleright \bullet" marking]
    & y
  \end{tikzcd}
\end{equation*}

Our proof will rest on the factorization $\sd(\Lambda^{2}_{0}) \hookrightarrow \sd(\Delta^{2})$ pictured in \hyperref[fig:factorization]{Figure~\ref*{fig:factorization}}. Denote the underlying factorization of simplicial sets by
\begin{equation*}
  A_{1} \overset{i_{1}}{\hookrightarrow}
  A_{2} \overset{i_{2}}{\hookrightarrow}
  A_{3} \overset{i_{3}}{\hookrightarrow}
  A_{4} \overset{i_{4}}{\hookrightarrow}
  A_{5} \overset{i_{5}}{\hookrightarrow}
  A_{6}
\end{equation*}

For each of the $A_{i}$ above, denote by
\begin{equation}
  \label{eq:inclusion_of_functors_respecting_labelling}
  \Fun'(A_{i}, \category{C}) \subseteq \Fun(A_{i}, \category{C})
\end{equation}
the full simplicial subset on those functors $A_{i} \to \category{C}$ which respect each labelling in \hyperref[fig:factorization]{Figure~\ref*{fig:factorization}}:
\begin{itemize}
  \item ingressive
  \item egressive
  \item $p$-cartesian
  \item $p$-cocartesian
  \item $p^{\dagger}$-cartesian
  \item $p_{\dagger}$-cocartesian
\end{itemize}
Note that because equivalences in $\category{C}$ belong to each of these classes of morphisms, the inclusion in \hyperref[eq:inclusion_of_functors_respecting_labelling]{Equation~\ref*{eq:inclusion_of_functors_respecting_labelling}} is an inclusion of connected components.

Similarly, we will denote by $\Fun'(A_{i}, \category{D})$ the full simplicial subset on functors which respect the ingressive and egressive labellings.

\begin{lemma}
  \label{lemma:homotopy_pullback_for_filling_procedure}
  The square
  \begin{equation*}
    \begin{tikzcd}
      \Fun'(A_{6}, \category{C})
      \arrow[r]
      \arrow[d]
      & \Fun'(A_{1}, \category{C})
      \arrow[d]
      \\
      \Fun'(A_{6}, \category{D})
      \arrow[r]
      & \Fun'(A_{1}, \category{D})
    \end{tikzcd}
  \end{equation*}
  is homotopy pullback.
\end{lemma}
\begin{proof}
  The factorization of \hyperref[fig:factorization]{Figure~\ref*{fig:factorization}} gives us a factorization of the above square into five squares
  \begin{equation*}
    \begin{tikzcd}
      \Fun'(A_{k+1}, \category{C})
      \arrow[r]
      \arrow[d]
      & \Fun'(A_{k}, \category{C})
      \arrow[d]
      \\
      \Fun'(A_{k+1}, \category{D})
      \arrow[r]
      & \Fun'(A_{k}, \category{D})
    \end{tikzcd},
    \qquad1 \leq k < 6,
  \end{equation*}
  each corresponding to one of the inclusions $i_{k}$. We will be done if we can show that each of these squares is homotopy pullback. It suffices to show that for each $k$, the map
  \begin{equation*}
    j_{k}\colon \Fun'(A_{k+1}, \category{C})
    \to
    \Fun'(A_{k}, \category{C})
    \times_{\Fun'(A_{k+1}, \category{D})}
    \Fun'(A_{k}, \category{D})
  \end{equation*}
  is a weak equivalence. First, we show that each of these maps is a Kan fibration. To see this, consider the square
  \begin{equation*}
    \noindent\makebox[\textwidth]{%
      \begin{tikzcd}[ampersand replacement=\&]
        \Fun'(A_{k+1}, \category{C})
        \arrow[r, hook]
        \arrow[d, swap, "j_{k}"]
        \& \Fun(A_{k+1}, \category{C})
        \arrow[d]
        \\
        \Fun'(A_{k}, \category{C})
        \times_{\Fun'(A_{k+1}, \category{D})}
        \Fun'(A_{k}, \category{D})
        \arrow[r, hook]
        \& \Fun(A_{k}, \category{C})
        \times_{\Fun(A_{k+1}, \category{D})}
        \Fun(A_{k}, \category{D})
      \end{tikzcd}.
    }
  \end{equation*}
  The right-hand map is a Kan fibration because of \hyperref[cor:kan_fib_and_isofib]{Corollary~\ref*{cor:kan_fib_and_isofib}}, and we have already seen that the top map is an inclusion of connected componenents. Hence each $j_{k}$ is a Kan fibration by \hyperref[cor:connected_components_kan_fibration]{Corollary~\ref*{cor:connected_components_kan_fibration}}. Thus, in order to show that each $j_{k}$ is a weak equivalence, it suffices to show that the fibers are contractible.

  First consider the case $k = 1$. For any map $\alpha$ below, consider the following pullback square.
  \begin{equation*}
    \begin{tikzcd}
      F_{\alpha}
      \arrow[r]
      \arrow[d]
      & \Fun'(A_{1}, \category{C})
      \arrow[d, "j_{1}"]
      \\
      \Delta^{0}
      \arrow[r, "\alpha"]
      & \Fun'(A_{0}, \category{C})
      \times_{\Fun'(A_{1}, \category{D})}
      \Fun'(A_{1}, \category{D})
    \end{tikzcd}.
  \end{equation*}
  The fiber $F_{\alpha}$ over $\alpha$ is the space of ways of completing a diagram of shape $A_{0}$ in $\category{C}$ to a diagram of shape $A_{1}$ in $\category{C}$ given a diagram of shape $A_{1}$ in $\category{D}$. This is the space of ways of filling $\Lambda^{2}_{2} \hookrightarrow \Delta^{2}$ in $\category{C}^{\dagger}$ lying over a $2$-simplex in $\category{D}^{\dagger}$. This is contractible by \hyperref[lemma:cartesian_horn_fillings_contractible]{Lemma~\ref*{lemma:cartesian_horn_fillings_contractible}}.

  The case $k = 2$ is similar, using \hyperref[lemma:inner_fib_contractible_fibers]{Lemma~\ref*{lemma:inner_fib_contractible_fibers}}.

  The case $k = 3$ is similar, using \hyperref[lemma:twice_cartesian_lifts_contractible]{Lemma~\ref*{lemma:twice_cartesian_lifts_contractible}}.

  The cases $k = 4$ and $k = 5$ use the dual to \hyperref[lemma:cartesian_horn_fillings_contractible]{Lemma~\ref*{lemma:cartesian_horn_fillings_contractible}}.
\end{proof}
We are now ready to show that morphisms of the form \hyperref[eq:form_of_p_cocartesian_morphisms]{Diagram~\ref*{eq:form_of_p_cocartesian_morphisms}} are $\pi'$-cocartesian. 

\begin{proposition}
  For any morphism $e$ of the form given in \hyperref[eq:form_of_p_cocartesian_morphisms]{Equation~\ref*{eq:form_of_p_cocartesian_morphisms}}, the square
  \begin{equation*}
    \noindent\makebox[\textwidth]{%
      \begin{tikzcd}[ampersand replacement=\&]
        \Fun^{\aCart}(\sd(\Delta^{2}), \category{C}) \times_{\Fun^{\aCart}(\sd(\Delta^{\{0, 1\}}), \category{C})} \{e\}
        \arrow[r]
        \arrow[d]
        \& \Fun^{\aCart}(\sd(\Lambda^{2}_{0}), \category{C}) \times_{\Fun^{\aCart}(\sd(\Delta^{\{0, 1\}}), \category{C})} \{e\}
        \arrow[d]
        \\
        \Fun^{\aCart}(\sd(\Delta^{2}), \category{D}) \times_{\Fun^{\aCart}(\sd(\Delta^{\{0, 1\}}), \category{D})} \{\pi e\}
        \arrow[r]
        \& \Fun^{\aCart}(\sd(\Lambda^{2}_{0}), \category{D}) \times_{\Fun^{\aCart}(\sd(\Delta^{\{0, 1\}}), \category{D})} \{\pi e\}
      \end{tikzcd}
    }
  \end{equation*}
  is homotopy pullback.
\end{proposition}
\begin{proof}
  This square factors horizontally into the two squares
  \begin{equation*}
    \noindent\makebox[\textwidth]{%
      \begin{tikzcd}[ampersand replacement=\&, column sep=tiny]
        \Fun^{\aCart}(\sd(\Delta^{2}), \category{C}) \times_{\Fun(\sd(\Delta^{\{0, 1\}}), \category{C})} \{e\}
        \arrow[r]
        \arrow[d]
        \& \Fun'(A_{6}, \category{C}) \times_{\Fun(\sd(\Delta^{\{0, 1\}}), \category{C})} \{e\}
        \arrow[d]
        \\
        \Fun^{\aCart}(\sd(\Delta^{2}), \category{D}) \times_{\Fun(\sd(\Delta^{\{0, 1\}}), \category{D})} \{\pi e\}
        \arrow[r]
        \& \Fun'(A_{6}, \category{D}) \times_{\Fun(\sd(\Delta^{\{0, 1\}}), \category{C})} \{\pi e\}
      \end{tikzcd}
    }
  \end{equation*}
  and
  \begin{equation*}
    \noindent\makebox[\textwidth]{%
      \begin{tikzcd}[ampersand replacement=\&, column sep=tiny]
        \Fun'(A_{6}, \category{C}) \times_{\Fun(\sd(\Delta^{\{0, 1\}}), \category{C})} \{e\}
        \arrow[r]
        \arrow[d]
        \& \Fun'(A_{1}, \category{C}) \times_{\Fun(\sd(\Delta^{\{0, 1\}}), \category{C})} \{e\}
        \arrow[d]
        \\
        \Fun'(A_{6}, \category{D}) \times_{\Fun(\sd(\Delta^{\{0, 1\}}), \category{C})} \{\pi e\}
        \arrow[r]
        \& \Fun'(A_{1}, \category{D}) \times_{\Fun(\sd(\Delta^{\{0, 1\}}), \category{C})} \{\pi e\}
      \end{tikzcd}.
    }
  \end{equation*}

  The first is homotopy pullback because the bottom map is an inclusion of connected components, and the second condition of \hyperref[thm:main]{Theorem~\ref*{thm:main}} guarantees that the fiber over an ambigressive $\sd(\Delta^{2}) \to \category{D}$ belonging to $\Fun'(A_{6}, \category{D})$ whose restriction to $\sd(\Delta^{\{0, 1\}})$ is $\pi e$ is precisely an ambigressive Cartesian functor $\sd(\Delta^{2}) \to \category{C}$ whose restriction to $\sd(\Delta^{\{0, 1\}})$ is $e$.

  That the second is homotopy pullback follows immediately from \hyperref[lemma:homotopy_pullback_for_filling_procedure]{Lemma~\ref*{lemma:homotopy_pullback_for_filling_procedure}}.
\end{proof}

This proves \hyperref[thm:main]{Theorem~\ref*{thm:main}}.

%\subfile{application_to_presheaves.tex}

\appendix

\section{Appendix}

\subsection{A brief list of model structures}
\label{sss:a_brief_list_of_model_structures}

The category $\SSet$ carries two model structures of which we will make frequent use:
\begin{itemize}
  \item The \emph{Kan} model structure, which has the following description.
    \begin{itemize}
      \item The fibrations are the Kan fibrations.

      \item The cofibrations are the monomorphisms.

      \item The weak equivalences are weak homotopy equivalences.
    \end{itemize}

  \item The \emph{Joyal} model structure. We will not give a complete description, referring the reader to \cite[Sec.\ 2.2.5]{highertopostheory}. We will make use of the following properties.
    \begin{itemize}
      \item The cofibrations are monomorphisms.

      \item The fibrant objects are quasicategories, and the fibrations between fibrant objects are isofibrations, i.e.\ inner fibrations with lifts of equivalences (cf.\ \cite[Cor.\ 2.6.5]{highertopostheory}).
    \end{itemize}
\end{itemize}

The category $\D\op$ has a Reedy structure, which gives us, together with the Kan model structure, a model structure on the category $\Fun(\D\op, \SSet) = \SSSet$ with the following properties.
\begin{itemize}
  \item The cofibrations are monomorphisms.

  \item The weak equivalences are level-wise weak homotopy equivalences.

  \item The fibrations are \emph{Reedy fibrations} (\hyperref[def:reedy_fibration]{Definition~\ref*{def:reedy_fibration}}).
\end{itemize}

This model category has two important left Bousfeld localizations which we will need.

\begin{itemize}
  \item The \emph{Segal space model structure} on $\SSSet$ is the left Bousfeld localization of the Reedy model structure at the set of spine inclusions
    \begin{equation*}
      S = \left\{\left.\Delta^{\{0, 1\}} \amalg_{\Delta^{\{1\}}} \Delta^{\{1, 2\}} \amalg_{\Delta^{\{2\}}} \cdots \amalg_{\Delta^{\{n-1\}}} \Delta^{\{n-1, n\}} \hookrightarrow \Delta^{n}\ \right|\ n \geq 2\right\},
    \end{equation*}
    which has the following properties.
    \begin{itemize}
      \item The cofibrations are the monomorphisms.

      \item The fibrant objects are Segal spaces.

      \item Every Reedy weak equivalence is a weak equivalence in the Segal space model structure.
    \end{itemize}

  \item The \emph{complete Segal space model structure} is the left Bousfeld localization of the Reedy model structure at the set $S$ of spine inclusions together with the inclusion $\{0\} \hookrightarrow I$, where $I$ is the (nerve of the) walking isomorphism. The complete Segal space model structure has the following properties.
    \begin{itemize}
      \item The cofibrations are the monomorphisms.

      \item The fibrant objects are complete Segal spaces.

      \item Every weak equivalence in the Segal space model structure is a weak equivalence in the complete Segal space model structure.
    \end{itemize}
\end{itemize}

\begin{proposition}[Joyal--Tierney]
  \label{prop:joyal-tierney_quillen_equivalence}
  There is a Quillen equivalence 
  \begin{equation*}
    \begin{tikzcd}
      p_{1}^{*} : \SSet^{\mathrm{Joyal}} \longleftrightarrow \SSSet^{\mathrm{CSS}} : i_{1}^{*}
    \end{tikzcd}
  \end{equation*}
  between the Joyal model structure on simplicial sets and the complete Segal space model structure on bisimplicial spaces. Here, the functors $p_{1}^{*}$ and $i_{1}^{*}$ are defined as follows.
  \begin{itemize}
    \item The functor $p_{1}^{*}$ takes a simplicial set $X$ to the bisimplicial space $(p_{1}^{*}X)_{mn} = X_{m}$. That is, the bisimplicial set $p_{1}^{*}X$ is constant in the the vertical direction, and each row is equal to the simplicial set $X$.

    \item The functor $i_{1}^{*}$ takes a bisimplicial set $K$ to the simplicial set $(i_{1}^{*}K)_{n} = K_{n0}$. That is, in the language of \hyperref[sss:a_brief_list_of_model_structures]{Subsection~\ref*{sss:a_brief_list_of_model_structures}}, the simplicial set $i_{1}^{*}K$ is the `zeroth row' of $K$.
  \end{itemize}
\end{proposition}

\begin{note}
  In the notation of
\end{note}

For more information, the reader is directed to \cite{qcats_vs_segal_spaces}.

\subsection{Divisibility of bifunctors}
\label{sss:divisibility_of_bifunctors}

In this section, we recall some key results from \cite{qcats_vs_segal_spaces}. We refer readers there for more information.

Let $\odot\colon \mathcal{E}_{1} \times \category{E}_{2} \to \category{E}_{3}$ be a functor of $1$-categories. We will say that $\odot$ is \emph{divisible on the left} if for each $A \in \category{E}_{1}$, the functor $A \odot -$ admits a right adjoint $A \backslash -$. In this case, this construction turns out also to be functorial in $A$; that is, we get a functor
\begin{equation*}
  - \backslash -\colon \category{E}_{1}\op \times \category{E}_{3} \to \category{E}_{2}.
\end{equation*}

Analogously, $\odot\colon \mathcal{E}_{1} \times \mathcal{E}_{2} \to \mathcal{E}_{3}$ is \emph{divisible on the right} if for for each $B \in \category{E}_{2}$, the functor $- \odot B$ admits a right adjoint $- / B$. In this case we get a of two variables
\begin{equation*}
  - / -\colon \category{E}_{3} \times \category{E}_{2}\op \to \category{E}_{1}.
\end{equation*}

\begin{example}
  The reader may find it helpful to keep in mind the cartesian product 
  \begin{equation*}
    - \times -\colon \SSet \times \SSet \to \SSet. 
  \end{equation*}
  In this case, both $A \backslash X$ and $X / A$ are the mapping space $X^{A}$.
\end{example}

If $\odot$ is divisible on both sides, then there is a bijection between maps of the following types:
\begin{equation*}
  A \odot B \to X,\qquad A \to X / B,\qquad B \to A \backslash X.
\end{equation*}
In particular, this implies that the functors $X / -$ and $- \backslash X$ are mutually right adjoint.

If both $\category{E}_{1}$ and $\category{E}_{2}$ are finitely complete and $\category{E}_{3}$ is finitely cocomplete, then from a map $u\colon A \to A'$ in $\category{E}_{1}$, a map $v\colon B \to B'$ in $\category{E}_{2}$, and a map $f\colon X \to Y$ in $\category{E}_{3}$, we can build the following maps.
\begin{itemize}
  \item From the square
    \begin{equation*}
      \begin{tikzcd}
        A \odot B
        \arrow[r]
        \arrow[d]
        & A' \odot B
        \arrow[d]
        \\
        A \odot B'
        \arrow[r]
        & A' \odot B'
      \end{tikzcd}
    \end{equation*}
    we get a map
    \begin{equation*}
      u \odot' v\colon A \odot B' \amalg_{A \odot B} A' \odot B \to A' \odot B'.
    \end{equation*}

  \item From the square
    \begin{equation*}
      \begin{tikzcd}
        A' \backslash X
        \arrow[r]
        \arrow[d]
        & A \backslash X
        \arrow[d]
        \\
        A' \backslash Y
        \arrow[r]
        & A \backslash Y
      \end{tikzcd}
    \end{equation*}
    we get a map
    \begin{equation*}
      \langle u \backslash f \rangle\colon A' \backslash X \to A \backslash X \times_{A \backslash Y} A' \backslash Y
    \end{equation*}

  \item From the square
    \begin{equation*}
      \begin{tikzcd}
        X / B'
        \arrow[r]
        \arrow[d]
        & X / B
        \arrow[d]
        \\
        Y / B'
        \arrow[r]
        & Y / B
      \end{tikzcd}
    \end{equation*}
    we get a map
    \begin{equation*}
      \langle f / v \rangle\colon X / B' \to X / B \times_{Y / B} Y / B'.
    \end{equation*}
\end{itemize}

\begin{proposition}
  \label{prop:equivalent_lifting_problems}
  With the above notation, the following are equivalent adjoint lifting problems:
  \begin{equation*}
    \begin{tikzcd}
      A \odot B' \amalg_{A \odot B} A' \odot B
      \arrow[r]
      \arrow[d]
      & X
      \arrow[d]
      \\
      A' \odot B'
      \arrow[r]
      \arrow[ur, dashed]
      & Y
    \end{tikzcd}
    \qquad
    \begin{tikzcd}
      A
      \arrow[r]
      \arrow[d]
      & X / B'
      \arrow[d]
      \\
      A'
      \arrow[r]
      \arrow[ur, dashed]
      & X / B \times_{Y / B} Y / B'
    \end{tikzcd}
  \end{equation*}
  \begin{equation*}
    \begin{tikzcd}
      B
      \arrow[r]
      \arrow[d]
      & A' \backslash X
      \arrow[d]
      \\
      B'
      \arrow[r]
      \arrow[ur, dashed]
      & A \backslash X \times_{A \backslash Y} A' \backslash Y
    \end{tikzcd}
  \end{equation*}
\end{proposition}

%\subfile{extra.tex}

\printbibliography

\end{document}